\documentclass[12pt]{article}
\usepackage{amsmath}
\usepackage{amsthm}
\usepackage{amssymb}
\usepackage{amscd}
\usepackage{amsfonts}
\usepackage{amsbsy}
\usepackage{graphicx}
\usepackage[dvips]{psfrag}
\usepackage{setspace}
\usepackage{indentfirst}
\usepackage{rotating,graphics,psfrag,epsfig}
\usepackage{float}
\usepackage{hyperref}
\usepackage{booktabs}
\usepackage{xcolor}
\usepackage{mathrsfs}
\usepackage{overpic}
\usepackage{doi}
\usepackage{tikz}
\usepackage{lineno}
\usepackage{dsfont}

\newcommand{\mom}[2]{\mathcal{M}_{#1}[#2]}

\newtheorem{teo}{Theorem}[section]
\newtheorem{corol}[teo]{Corollary}

\newtheorem{propo}[teo]{Proposition}
\newtheorem{lem}[teo]{Lemma}
\newtheorem{defi}[teo]{Definition}

\newtheorem{obs}[teo]{Remark}

\newtheorem{main}{Theorem}

\renewenvironment{proof}{\noindent {\emph{Proof.}}}{\qed\vspace{0.4cm}}

\begin{document}

%\linenumbers
%\onehalfspacing

\begin{center}
{\Large{\bf On the multiplicity of crossing limit cycles of double diagonalizable piecewise linear vector fields}}
\end{center}

\bigskip

\bigskip

\begin{center}
{\large Lucas Queiroz Arakaki$^\ast$},\\
{\em IMECC - Universidade Estadual de Campinas,\\
Campinas, SP, Brazil.} \\
e-mail: larakaki@unicamp.br
\end{center}
\begin{center}
{\large Daniel Panazzolo},\\
{\em IRIMAS - Universit\'e de Haute-Alsace\\ 68093 Mulhouse, France.\\}
e-mail: daniel.panazzolo@uha.fr
\end{center}

\begin{center}
{\bf Abstract} 
\end{center}
We consider planar piecewise linear vector fields with two zones separated by a straight line and having two diagonalizable components, which we call the \emph{double diagonalizable configuration}. In this configuration, the Poincaré half-maps correspond to separatrices of saddle points in the foliations of symmetric differential forms $\omega_r(x)-\omega_r(y)$, where $\omega_r$ is a rational $1$-form with exactly three simple poles. Limit cycles correspond to non-diagonal intersections of two such separatrices, and the multiplicity of a limit cycle is precisely the order of contact at the intersection point. We prove that two distinct separatrices which are not nested segments of the same straight line have a contact of order at most $4$ at a diagonal point and at most $2$ at a non-diagonal point. As a consequence, we obtain that the multiplicity of limit cycles in the double diagonalizable configuration is at most two, this bound being sharp.
\bigskip

\noindent {\small {\bf Key-words}: Planar piecewise linear vector fields, Limit cycles, Upper bounds, Poincaré half-maps, Differential forms, Foliations}

\noindent {\small {\bf 2020 Mathematics Subject Classification:} 34A36, 34C07, 34C05, 37C27, 34A26}

\noindent \makebox [40mm]{\hrulefill}

\noindent {\footnotesize $^\ast$Corresponding author}.

\tableofcontents

\section{Introduction and statement of the main results}

The qualitative study of differential systems has been one of the main interests of the mathematical community for more than a century, which is a consequence of Hilbert's 16th Problem whose object is to determine, if it exists, the upper bound on the number of limit cycles of planar polynomial vector fields. Stemming from its original formulation, several different versions have been considered. One of these versions consists of studying, instead of polynomial vector fields, nonsmooth vector fields.

The investigation of limit cycles in nonsmooth planar vector fields itself has a history that dates back to the 1930s with the work of Andronov et al. \cite{Andronov37}. Perhaps one of the most studied subclasses of nonsmooth vector fields is that of the planar piecewise linear vector fields, with two zones separated by a straight line. In other words, up to an affine change of variables, vector fields of the form
%\begin{equation}
%\dot{\mathbf{x}}=\left\{\begin{array}{l} A_{-}\mathbf{x}+\mathbf{b}_{-},\; y< 0,\\
%A_{+}\mathbf{x}+\mathbf{b}_{+},\; y> 0,
%\end{array}\right.
%\end{equation}
\begin{equation}\label{eq:1}
X=\mathds{1}_{y>0} X_+ + \mathds{1}_{y<0} X_-,
\end{equation}
where $X_{\pm}$ are linear vector fields, and $\mathds{1}_{\mathcal{S}}$ is the characteristic function of a set $\mathcal{S}$. The dynamics on the \emph{switching set} $\Sigma=\{y=0\}$ are defined by the Filippov convention \cite{Filippov}. In this setting, a \emph{limit cycle} (also referred to as \emph{crossing limit cycle}) is a crossing non-trivial closed orbit which is isolated from other closed orbits.

Several authors have studied system \eqref{eq:1} in the last decades and significant progress has been made concerning the upper and lower bounds on the number of limit cycles. Freire et al. in \cite{Freire98} proved that \eqref{eq:1} has at most one limit cycle in the case that the Filippov vector field is continuous, which was first conjectured by Lum and Chua \cite{LumChua}. In the discontinuous case, several upper bounds have been obtained for particular cases of system \eqref{eq:1}, for instance: in \cite{LliNovMAT15}, it was proved that the maximum number of limit cycles is two if one of the components has a focus at $\Sigma$, and one if either component has a center or a weak saddle (a saddle with hyperbolicity ratio 1). The upper bound is also one for refracting systems \cite{LiLiuLli} and for sewing systems, that is, systems without sliding region \cite{MedradoTorre15,CarSanNovUnique}.

Recently, Carmona, Fernández-Sánchez and Novaes \cite{CarSanNov23} established that the uniform upper bound on the number of limit cycles of \eqref{eq:1} is $L^\ast\leqslant 8$. Their breakthrough was based on the integral characterization of Poincaré half-maps \cite{CarmonaIntegral} and an extension of Khovanskii's theory, and was further summarized into a methodology to determine upper bounds on limit cycles of piecewise linear vector fields in \cite{CarSanNovpreprint}. However, the best known lower bound on the maximum number of limit cycles of \eqref{eq:1} is three (see, for instance, \cite{HuanYang12,LliPonce12,BuzziPessoa,FreireMechanism14}). The known examples of piecewise linear vector fields with three limit cycles involve at least one focus, whereas for configurations without foci the best lower bound is two \cite{ChenNodes,LliMATJoan}.

Let us now formulate our main result. We say that a piecewise linear vector field $X$ of the form \eqref{eq:1} is \emph{double diagonalizable} (or of \emph{D-D type}, for short) if both its components $X_+$ and $X_-$ have a diagonalizable linear part with non-zero eigenvalues. Equivalently, each component $X_{\pm}$ of $X$ has a saddle or a proper node singularity, which can be either visible or invisible. The particular case of D-D piecewise linear vector fields where the components have facing saddles has been considered in \cite{ArtesMAT,LliValls21,ShaoRefractive}. For two nodes, Chen and Liu \cite{ChenNodes} proved, by a computer-assisted analysis, that there are at most two limit cycles counting multiplicities. Other types of D-D vector fields were considered in \cite{LliMATJoan}.

%Our objective in this paper is to determine the maximum multiplicity of any limit cycle in the double diagonalizable configuration. More precisely, our main result is the following:

%\begin{main}\label{Teo:A}
%There are no limit cycles of multiplicity 3 for piecewise linear vector fields with two diagonalizable components.
%\end{main}
The main result of the paper is the following Theorem:

\begin{main}\label{Teo:A}
All limit cycles of piecewise linear vector fields of D-D type have multiplicity at most two. Moreover, this bound is sharp.
\end{main}

To our knowledge, Theorem \ref{Teo:A} is new for the saddle--saddle and node--saddle configurations whose fold points do not coincide. Indeed, in the node--node case it follows from \cite{ChenNodes}, and when the two fold points coincide there is no sliding region and every limit cycle is hyperbolic by \cite{CarSanNovUnique}. Our proof is analytic and treats all D-D configurations at once. Since the displacement function depends analytically on the parameters, Theorem \ref{Teo:A} also implies that at most two limit cycles bifurcate from a crossing limit cycle of a D-D vector field under small perturbations within the class \eqref{eq:1}.

%Our approach is based on the study of the contacts of separatrices of foliations of triple-polar differential forms. In brief, due to the diagonalizability hypothesis, for each component $X_{\pm}$ a Poincaré half-map $x\mapsto \sigma_\pm(x)$ is well-defined on some interval contained in $\Sigma$. Moreover, each component $X_\pm$ has a first integral of Darboux type. Through this first integral, one can express the graphs of $\sigma_\pm$  as separatrices associated to  symmetric differential forms of triple polar type. More precisely, the curves $(x,\sigma_\pm(x))$ can be interpreted as one of the separatrices of a saddle point (see Figure \ref{Fig:RF-triplepolar}) in the foliation of a symmetric differential form $\Omega_r=\omega_r(x)-\omega_r(y)$, where $\omega_r$ is a \emph{triple polar} $1$-form, in simpler terms, a differential form given by
%\[\omega_r=\frac{z-r_0}{(z-r_1)(z-r_2)}dz,\]
%where $r=(r_0,r_1,r_2)$ is a triple of distinct real numbers. We were able to prove that, for given triples $r$ and $r'$, if the correspondent separatrices $\Gamma_r$ and $\Gamma_{r'}$ have a contact of order 5 at a point of the diagonal $\{x=y\}$ or have a contact of order 3 at a non-diagonal point, then they must either coincide or be nested segments of the same straight line.

Due to the diagonalizability hypothesis, for each component $X_{\pm}$ a Poincaré half-map $x\mapsto \sigma_\pm(x)$ is defined on some interval contained in $\Sigma$, and $X_\pm$ has a first integral of Darboux type. Through this first integral, one can express the graphs of these maps as separatrices associated with suitably defined differential forms. More precisely, parametrizing the discontinuity line $\Sigma$ by the $x$-coordinate and denoting by $\Sigma_{\rm in}^\pm$ and $\Sigma_{\rm out}^\pm$ the domain and image of the Poincaré half-map of $X_\pm$ respectively, the curves $\{(x,y) \mid x \in \Sigma_{\rm in}^\pm, y=\sigma_\pm(x)\}$ can be interpreted as separatrices of a saddle point (see Figure \ref{Fig:RF-triplepolar}) within the foliation of the symmetric differential form $\Omega_r=\omega_r(x)-\omega_r(y)$. Here,$$\omega_r(z) = \frac{z-r_0}{(z-r_1)(z-r_2)}dz,$$and $r=(r_0, r_1, r_2)$ is a triple of distinct real numbers. Using this correspondence, each limit cycle of multiplicity $n$ for a discontinuous vector field $X$ of the form \eqref{eq:1} corresponds precisely to a point of contact of order $n$ between the separatrices of the corresponding symmetric forms.
%We remark that Theorem \ref{Teo:A} encompasses piecewise linear vector fields whose components have proper nodes and saddles. Its bound on the multiplicity for the particular case where both components have proper nodes was previously obtained by Chen and Liu \cite{ChenNodes}, using a computer assisted method.

The present work is structured as follows: In Section \ref{S:Half-return}, we establish the above-mentioned connection between  the Poincaré half-maps of piecewise linear vector fields of D-D type  and separatrices associated to symmetric differential forms $\Omega_r$ as defined above. From this point, we state the intermediate result (Theorem \ref{Teo:I}) that if the order of contact of two such separatrices is either $5$ at a diagonal point or $3$ at a non-diagonal point, then these separatrices either coincide or are nested line segments. In Section \ref{S:contact}, we prove Theorem \ref{Teo:I} by studying the contact at diagonal points in subsection \ref{S:5contact}, whereas subsection \ref{S:3contact} is dedicated to the investigation of $3$-contacts at non-diagonal points which is translated via a projective change of coordinates to the determination of the fibers of an analytic map $\Theta:\mathbf{D}\subset\mathbb{R}^2\to\mathbb{R}^2$, where $\mathbf{D}$ is an appropriate open subset of $\mathbb{R}^2$. Some illustrative examples are presented in Section \ref{S:Examples}. Section \ref{S:Proofs} proves all technical results needed to describe the fibers of the map $\Theta$. As we shall see, this description is fully determined by the sign of the Jacobian determinant $\det D\Theta$. To determine this sign, we express all relevant functions as parameter-dependent weighted integrals over the real line. This allows us to reduce the evaluation of $\operatorname{sign}(\det D\Theta)$ to comparing the growth rates of two moment ratios associated with these integral representations.

\section{Half-return maps via Darboux first integrals}\label{S:Half-return}

We start by discussing the general concept of half-return maps for smooth vector fields defined on domains with boundary. 
Let $X$ be a smooth vector field defined on an open subset $U\subset\mathbb{R}^n$. Given a closed subset $\mathcal{V}\subset U$ with boundary $\Sigma=\partial\mathcal{V}$, we define the \emph{domain of the half-return map of $X\vert_\mathcal{V}$ at $\Sigma$} as the set
\[\Sigma_{\rm in}=\{p\in\Sigma: \exists\,t_p>0: \varphi_X(t_p;p)\in\Sigma\,\text{and }\varphi_X((0,t_p);p)\subset\mathring{\mathcal{V}}\},
\]
where $\varphi_X(\cdot;p)$ denotes the flow of $X$ through the point $p$ with $\varphi_X(0;p)=p$ and $\mathring{\mathcal{V}}$ denotes the interior of the set $\mathcal{V}$. In this case, the map $\sigma:\Sigma_{\rm in}\to\Sigma$ given by $\sigma(p)=\varphi_X(t_p,p)$ is called the \emph{half-return map} (or Poincaré half-map) of $X\vert_\mathcal{V}$ at $\Sigma$. We denote by $\Sigma_{\rm out}=\sigma(\Sigma_{\rm in})$ its image. Whenever $\Sigma_{\rm in}$ is non-empty, we say that a half-return map is \emph{well-defined} for $X\vert_\mathcal{V}$. 

\subsection{Half-return maps for linear vector fields in half-planes}

We want to completely characterize the domains $\Sigma_{\rm in}$ and images $\Sigma_{\rm out}$ in the case that $X$ is a diagonalizable linear vector field and $\mathcal{V}$ is a closed half-plane. 

We say that a planar linear vector field $X$ is \emph{diagonalizable} if its Jordan canonical form is diagonal. In other words, there exists an affine coordinate change that transforms $X$ into $X=\lambda_1u\,\partial_u+\lambda_2 v\,\partial_v$,
with $\lambda_1\lambda_2\neq 0$. Up to an additional time rescaling $t\to t/\lambda_1$, we see that any diagonalizable vector field is orbitally affinely equivalent to 
\begin{equation}\label{eq:Jordan}
X_J=u\partial_u-\lambda v\partial_v,
\end{equation}
and  the number $\lambda \in \mathbb{R}^\ast$ is called the hyperbolicity ratio of $X$.

In this setting, the half-return maps, also called {\em Poincaré half-maps}, have been extensively studied in the literature, mainly connected to the dynamics of piecewise smooth vector fields (see \cite{CarmonaIntegral,CarSanNov23,Freire98,Freire2012}). A recent result, due to Carmona and Fern{\'a}ndez-S{\'a}nchez \cite{CarmonaIntegral}, characterizes these maps as the Cauchy Principal Value of the integral of a rational function. Their construction is based on a Li{\'e}nard canonical form. In Subsection \ref{SS:FI}, we will take a different approach, constructing the half-return map from the Darboux first integrals of diagonalizable linear vector fields. In order to do so, we need some preliminary results. 

Let us initially consider the following simplified situation. Let $X$ be given in the orbital Jordan form \eqref{eq:Jordan} for some $\lambda\in\mathbb{R}\setminus[-1,0]$ and $\Sigma$ be given by the affine line $f(u,v)=0$, where $f(u,v)=u+v-1$. Then, consider the lines
\[\ell_1=\{u=0\},\quad \ell_2=\{v=0\},\quad Xf(u,v)\equiv 0,\]
where $Xf$ denotes the Lie derivative of $f$ with respect to $X$. We consider also the intersection points $p=\left(\frac{\lambda}{\lambda+1},\frac{1}{\lambda+1}\right)=\Sigma\cap\{Xf\equiv0\}$, $p_1=(0,1)=\Sigma\cap\ell_1$ and $p_2=(1,0)=\Sigma\cap\ell_2$.

The following result characterizes the domain and image of the half-return maps of $X$ restricted to the half-planes $\Sigma^+=\{f(u,v)\geqslant 0\}$ and $\Sigma^-=\{f(u,v)\leqslant 0\}$. In the statement, for any two points $q_0, q_1 \in \Sigma$, $(q_0, q_1)$ denotes the open segment of $\Sigma$ between $q_0$ and $q_1$. The notation $(q_0,\infty)$ indicates the infinite ray starting from $q_0$.

\begin{lem}\label{Lema:half-return-Jordan}
In the above setting, the following holds:
\begin{enumerate}
    \item If $\lambda>0$ (saddle case):
    \begin{itemize}
        \item[a)] If $\mathcal{V}=\Sigma^+$, then $\Sigma_{\rm in}=\emptyset$;
        \item[b)] If $\mathcal{V}=\Sigma^-$, then $\Sigma_{\rm in}=(p_1,p)$ and $\Sigma_{\rm out}=(p,p_2)$.
    \end{itemize}
    \item If $\lambda<-1$ (node case):
    \begin{itemize}
        \item[a)] If $\mathcal{V}=\Sigma^-$, then $\Sigma_{\rm in}=\emptyset$;
        \item[b)] If $\mathcal{V}=\Sigma^+$, then $\Sigma_{\rm in}=(p_2,p)$ and $\Sigma_{\rm out}=(p,\infty)$
    \end{itemize}
\end{enumerate}
\end{lem}

\begin{figure}[!ht]
    \centering
    \vspace*{0.5cm}
    \begin{overpic}[width=0.9\linewidth]{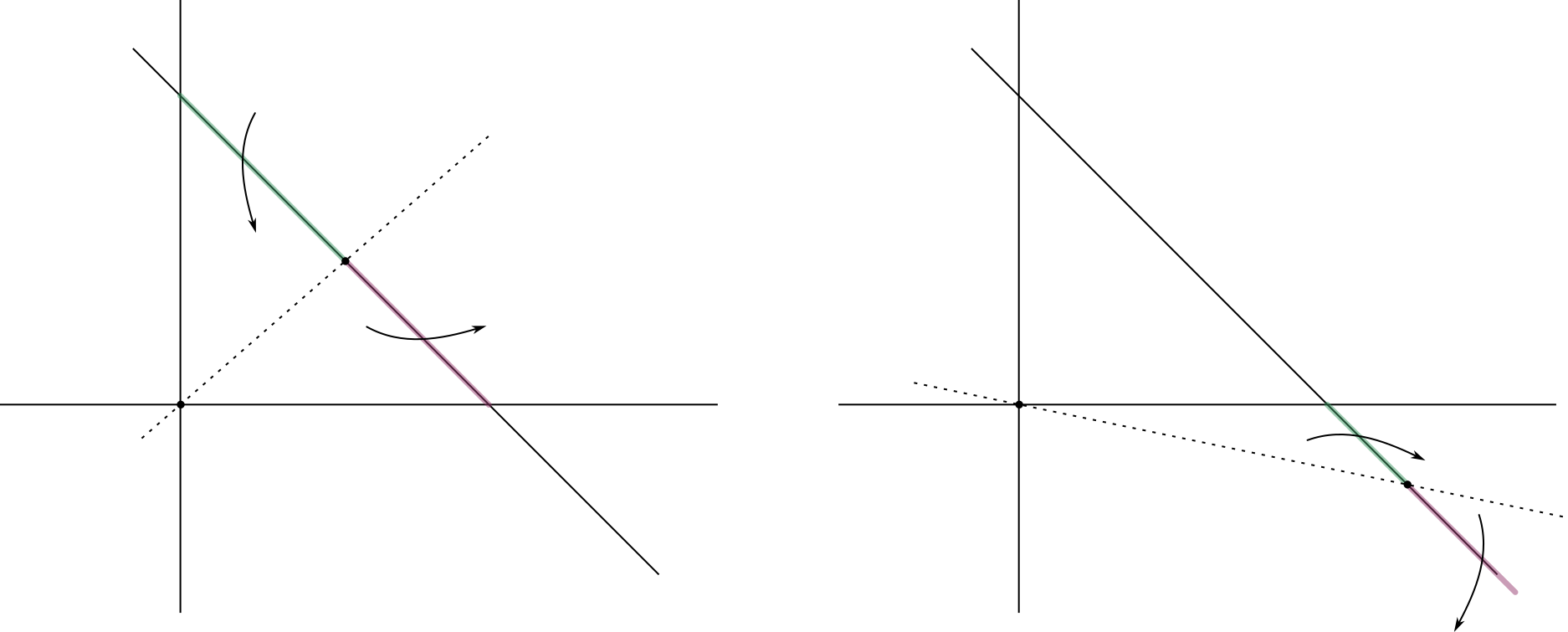}
     % --- highlights (coordinates in overpic units: image width = 100) ---
    \put(0,0){\begin{tikzpicture}[overlay, x=\unitlength, y=\unitlength,
                                  line width=4pt, line cap=round, opacity=.5]
        % Panel 1
        \draw[green!60!black]   (11.6,34.0) -- (22.1,23.4); % Sigma_in
        \draw[magenta!70!black] (22.1,23.4) -- (31.1,14.3); % Sigma_out
        % Panel 2
        \draw[green!60!black]   (84.5,14.3) -- (89.7, 9.2); % Sigma_in
        \draw[magenta!70!black] (89.7, 9.2) -- (96.0, 2.6); % Sigma_out
    \end{tikzpicture}}
    % --- labels ---
    \put(12.6,34.7){$p_1$}
    \put(28.4,11.4){$p_2$}
    \put(65.6,34.7){$p_1$}
    \put(85.6,15.4){$p_2$}
    \put(23,40){1. $\lambda>0$:}
    \put(75,40){2. $\lambda<-1$:}
    \put(23.5,23.5){$p$}
    \put(87.5,7.5){$p$}
    \put(20,27.5){$\Sigma_{\rm in}$}
    \put(32,16){$\Sigma_{\rm out}$}
    \put(90,12){$\Sigma_{\rm in}$}
    \put(96,4){$\Sigma_{\rm out}$}
    \end{overpic}
    \caption{The domain $\Sigma_{\rm in}$ and image $\Sigma_{\rm out}$ of the half-return map in the setting of Lemma \ref{Lema:half-return-Jordan}}
    \label{Fig:Lemma-Jordan}
\end{figure}

The proof of this lemma is immediate from a standard analysis of linear vector fields and an application of barrier functions. Figure \ref{Fig:Lemma-Jordan} illustrates the simplified situation. We now take on the more general setting:
\begin{enumerate}
    \item $X$ is a diagonalizable linear vector field;% with distinct hyperbolicity ratio $\lambda\neq -1$;
    \item $\mathcal{V}=\{f\geqslant0\}$ is an affine half-plane with boundary $\Sigma=\{f=0\}$, where $f$ is a non-zero affine function.
\end{enumerate}
%Condition 1 above excludes star nodes, since those do not admit half-return maps. 

\begin{figure}[!ht]
    \centering
    \begin{overpic}[width=0.4\linewidth]{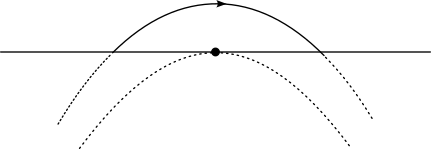}
    \put(50,17.5){$p$}
    \end{overpic}
    \caption{An invisible fold.}
    \label{Fig:fold}
\end{figure}

We say that $p\in\Sigma$ is an \emph{invisible fold} of $X\vert_\mathcal{V}$ if the lines $\Sigma$ and $Xf\equiv0$ intersect only at the point $p$ and additionally $X^2f(p)<0$ (Figure \ref{Fig:fold}).

\begin{propo}\label{Propo:half-return}
In the above setting, $X\vert_{\mathcal{V}}$ has a well-defined half-return map (i.e.~$\Sigma_{\rm in}\ne \emptyset$) if and only if $X\vert_\mathcal{V}$ has an invisible fold. Moreover, in this case, up to an affine change of coordinates and time rescaling, we are precisely in one of the cases (1.b) or (2.b) listed in Lemma \ref{Lema:half-return-Jordan}, with the roles of $\Sigma_{\rm in}$ and $\Sigma_{\rm out}$ interchanged when such rescaling reverses time. In particular, the fold point is precisely the point with coordinates $ (\frac{\lambda}{\lambda+1},\frac{1}{\lambda+1})$.
\end{propo}
\begin{proof}
First, we observe that if the lines $Xf\equiv 0$ and $\Sigma$ coincide, then $\Sigma$ is necessarily an invariant line of $X$ and no half-return map can be defined. Similarly, if these lines do not intersect, then $Xf$ has a fixed sign on $\Sigma$ and no half-return map for $X\vert_\mathcal{V}$ can be defined.

Analogously, if $Xf\equiv 0$ and $\Sigma$ intersect at the equilibrium of $X$ then the invariant lines of $X$ create a barrier between the orbits starting at $\Sigma\cap\{Xf>0\}$ and the set $\Sigma\cap\{Xf<0\}$. Therefore, $\Sigma_{\rm in}=\emptyset$ in this case, as well.

Thus, up to an affine change of coordinates and time rescaling, we can assume that $X$ is in the orbital Jordan form \eqref{eq:Jordan} and $f(u,v)=u+v-1$ or $f(u,v)=1-u-v$. In other words, we are precisely in the setting of Lemma \ref{Lema:half-return-Jordan}, with either $\mathcal{V}=\Sigma^+$ or $\mathcal{V}=\Sigma^-$.

To conclude, it suffices to remark that, in the setting of Lemma \ref{Lema:half-return-Jordan}, the half-return map of $X\vert_\mathcal{V}$ has a non-empty domain $\Sigma_{\rm in}$ if and only if $X\vert_\mathcal{V}$ has an invisible fold at the tangency point $p$.
\end{proof}

\begin{obs}
In the case where $\Sigma =\{y=0\}$ and $X$ assumes the general form \[X=(a_{11}x+a_{12}y+b_1)\partial_x+(a_{21}x+a_{22}y+b_2)\partial_y,\]  the restricted vector field  $X\vert_{y>0}$ has an invisible fold point at $\Sigma$ if and only if 
\begin{equation*}%\label{eq:invisible fold}
a_{21}\neq 0,\quad (a_{21}b_1-a_{11}b_2)<0.
\end{equation*}
By Proposition \ref{Propo:half-return}, when $X$ is diagonalizable these are necessary and sufficient conditions for $X$ restricted to $\{y\geqslant 0\}$ to have a well-defined half-return map on $\Sigma$. 
\end{obs}

\subsection{Separatrices of triple polar forms}\label{SS:Separatrices}

In this subsection, we study the intersection of separatrices of triple polar forms and their order of contact. As we shall see in Subsection \ref{SS:FI}, the half-return maps of each component of a piecewise linear D-D vector field are given precisely by separatrices of this kind. Theorem \ref{Teo:A} will follow directly from the fact that, unless two distinct separatrices are nested segments of the same straight line, their order of contact is at most four at a point of the diagonal $\{x=y\}$ and at most two at a non-diagonal point. We now formalize these ideas.

\begin{defi}[Triple polar forms]
We say that a rational 1-form $\omega$ on $\mathbb{CP}^1$ is a \emph{triple polar form} when it has precisely three simple poles. 
\end{defi}

\begin{figure}[!ht]
    \centering
    \vspace*{0.5cm}
    \begin{overpic}[width=\linewidth]{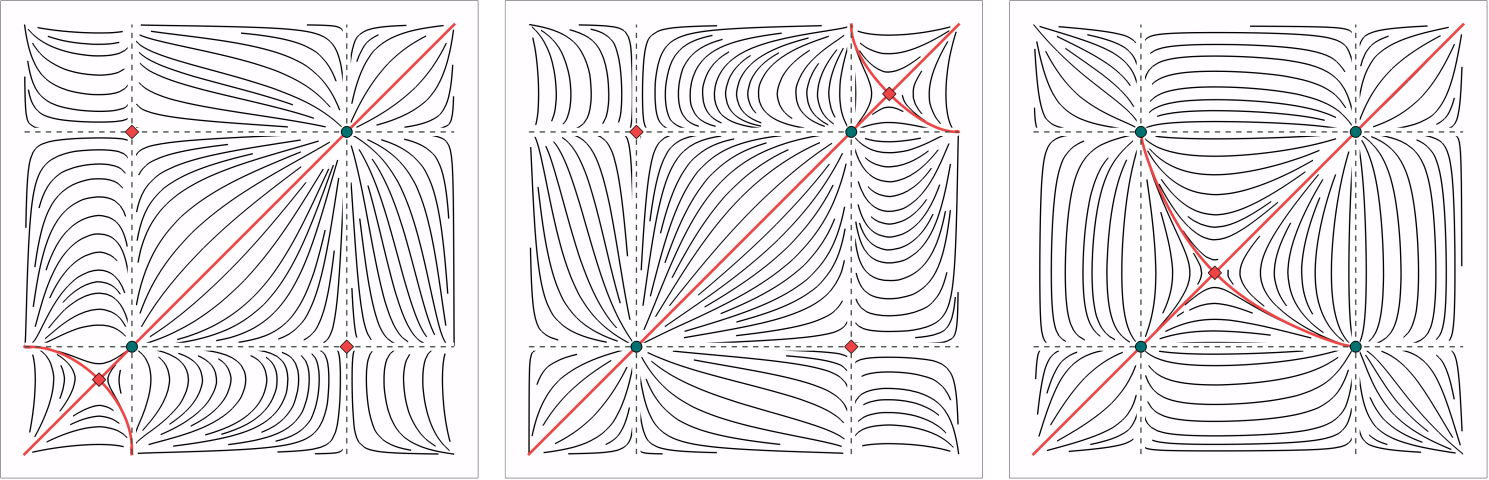}
    \put(10,33.5){(a) $r_0<r_1$}
    %\put(3.5,-1.5){$r_0$}
    \put(8,-1.5){$r_1$}
    \put(23,-1.5){$r_2$}
    \put(43,33.5){(b) $r_0>r_2$}
    %\put(62,-1.5){$r_0$}
    \put(42,-1.5){$r_1$}
    \put(57,-1.5){$r_2$}
    \put(75,33.5){(c) $r_1<r_0<r_2$}
    %\put(81.5,-1.5){$r_0$}
    \put(76,-1.5){$r_1$}
    \put(91,-1.5){$r_2$}
    \end{overpic}
    \caption{The foliation of $\Omega_r$ for different configurations of $r_0,r_1,r_2$ in $\mathbb{M}$. The manifold $\mathbb{M}$ is diffeomorphic to a 2-torus, which is topologically obtained by gluing the opposite sides of these squares. The separatrices of the saddle $\mathbf{s}$ are highlighted as thick red lines. Colors available online.}
    \label{Fig:RF-triplepolar}
\end{figure}

Any triple polar form can be transformed via a M\"obius change of variables and multiplication by a nonzero constant into the following normal form:
\begin{equation*}%\label{eq:triplepolar}
\omega_r(z)=\frac{z-r_0}{(z-r_1)(z-r_2)}dz,
\end{equation*}
where $r=(r_0,r_1,r_2)$ has distinct entries. In this normal form, the poles are located at $z=r_1$, $z=r_2$ and $z=\infty$. We say that $\omega_r$ is real if $r_0,r_1,r_2\in\mathbb{R}$. From now on, we will always assume that $\omega_r$ is real and $r_1<r_2$. Let us set up the following notations: let $\mathbb{M}=\mathbb{RP}^1\times\mathbb{RP}^1$, and $\Lambda=\{(r_0,r_1,r_2)\in\mathbb{R}^3:r_0\not\in \{r_1,r_2\},r_1<r_2\}$. The family $\mathbf{\Omega}=\{\Omega_r\}_{r\in\Lambda}$ of real symmetric 1-forms on $\mathbb{M}$ given by
\begin{equation}\label{eq:Omega_r}
\Omega_r=\omega_r(x)-\omega_r(y),
\end{equation}
will be called the \emph{real triple polar family}. The non-oriented analytic foliation associated to $\Omega_r$ on the affine plane $\mathbb{R}^2 \subset \mathbb{M}$ is symmetric with respect to the diagonal, and coincides with the (non-oriented) phase portrait of the differential system
\begin{equation}\label{eq:DE3}
\dot{x}=\frac{y-r_0}{(y-r_1)(y-r_2)},\quad \dot{y}=\frac{x-r_0}{(x-r_1)(x-r_2)}.
\end{equation}
The points $\{r_1,r_2,\infty\}$ delimit three distinct segments of $\mathbb{RP}^1$. Let $\mathscr{I}_r$ be the segment containing $r_0$. More precisely, 
\begin{equation}\label{eq:Ir}
\mathscr{I}_r=\left\{\begin{array}{l}
(\infty,r_1),\quad\text{for }r_0<r_1;\\
(r_2,\infty),\quad\text{for }r_0>r_2;\\
(r_1,r_2),\quad\text{for }r_1<r_0<r_2.
\end{array}\right.
\end{equation}
An elementary analysis of system \eqref{eq:DE3} concludes that inside the \emph{fundamental square} $\mathscr{I}_r^2$, system \eqref{eq:DE3} has a hyperbolic saddle at $\mathbf{s}=(r_0,r_0)$ with orthogonal associated eigenspaces. In particular, the saddle $\mathbf{s}$ has the diagonal $\{x=y\}$ as a separatrix and its second separatrix, henceforth denoted by $\Gamma_r$, is tangent to the line $x+y=2r_0$ at $\mathbf{s}$ (see Figure \ref{Fig:RF-triplepolar}).

We denote by $\mathbf{\Gamma}=\{\Gamma_r\}_{r\in\Lambda}$, the corresponding family of non-diagonal separatrices associated to $\mathbf{\Omega} = \{\Omega_r\}_{r\in\Lambda}$.

The following classical results and notions will allow us to characterize the separatrices $\Gamma_r\in\mathbf{\Gamma}$.

\begin{defi}
We say that a differentiable function $f:\mathscr{J}\to\mathbb{R}$ defined on an open interval $\mathscr{J}\subset\mathbb{R}$ is \emph{unimodal}, when there exists a unique critical point $x_c\in\mathscr{J}$ which is a global extremum point of $f$ on $\mathscr{J}$.
\end{defi}

\begin{figure}[!ht]
    \centering
    \begin{overpic}[width=0.5\linewidth]{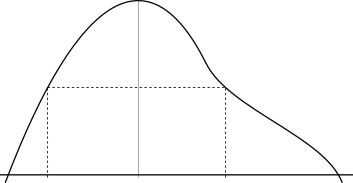}
    \put(40,30){$f(x)$}
    \put(3.5,-1.5){$x_0$}
    \put(92,-1.5){$x_1$}
    \put(39,-1.5){$x_c$}
    \put(12.7,-1.5){$x$}
    \put(63,-1.5){$\sigma(x)$}
    \end{overpic}
    \caption{The involution $\sigma$ defined by a unimodal function $f$.}
    \label{Fig:involution}
\end{figure}

\begin{lem}\label{Lema:unimodal-involution}
If $f:(x_0,x_1)\to\mathbb{R}$ is unimodal, and $\lim\limits_{x\to x_0^+}f(x)=\lim\limits_{x\to x_1^-}f(x)$, then there exists a unique continuous non-trivial involution $\sigma:(x_0,x_1)\to (x_0,x_1)$ such that $f\circ\sigma(x)=f(x)$. Moreover, 
\begin{itemize}
    \item[i.] Let $x_c$ be the unique critical point of $f$. Then the involution exchanges $(x_0,x_c)$ with $(x_c,x_1)$, i.e. the restriction $\sigma:(x_0,x_c)\to(x_c,x_1)$ is a bijection;
    \item[ii.]If, additionally, $f$ is analytic, then so is $\sigma$.
\end{itemize}
\end{lem}
\begin{proof}
Assume, without loss of generality, that $x_c$ is a global minimum of $f$. Let $L=\lim\limits_{x\to x_0^+} f(x)$, with $L\in\mathbb{R}\cup\{+\infty\}$. Each value $y\in(f(x_c),L)$ is attained exactly twice: once in $(x_0,x_c)$ and once in $(x_c,x_1)$. Define $\sigma(x_c)=x_c$ and $\sigma(x)=\bar{x}$, where $x\neq\bar{x}$ and $f(x)=f(\bar{x})$. Monotonicity of $f$ in each side of $x_c$ implies the uniqueness of $\sigma$. Moreover, $\sigma\circ\sigma(x)=\sigma(\bar{x})=x$ and $\sigma$ is indeed an involution and exchanges $(x_0,x_c)$ with $(x_c,x_1)$ by definition, which proves assertion (i).

If $f$ is analytic, there exists $k\in\mathbb{N}$ such that
\[f(x)=f(x_c)+(x-x_c)^{2k}u(x),\quad u(x_c)>0.\]
Let $g(x)=(x-x_c)u(x)^{1/2k}$. Then, $\bar{\sigma}(x)=g^{-1}(-g(x))$ is a nontrivial involution satisfying $f\circ\bar\sigma(x)=f(x)$. By uniqueness, $\sigma=\bar{\sigma}$ and thus $\sigma$ is analytic near $x_c$. Analyticity on $(x_0,x_1)\setminus\{x_c\}$ is implied by the Implicit Function Theorem. The proof for the case in which $x_c$ is a global maximum is immediate by considering $-f$. This completes the proof of assertion (ii). Figure \ref{Fig:involution} illustrates the proof.
\end{proof}

In the remainder of our paper, we will refer to non-trivial involutions simply as \emph{involutions}. We recall that any smooth involution $\sigma$ defined on a connected interval is strictly decreasing and must satisfy $\sigma'(x_c)=-1$ at its fixed point $x_c$ (see, for instance \cite{ShishaInvolutions,ZampieriInvolution}). For a more detailed treatment of involutions and their relationship with differential systems, we refer the reader to \cite{CimaManVil}.\medskip

Given a parameter $r=(r_0,r_1,r_2) \in \Lambda$, let us now consider the function 
\begin{equation}\label{eq:A,B,h}
h_r(u):=A\log|u-r_1|+B\log|u-r_2|,
\end{equation}
where $A=\frac{r_1-r_0}{r_1-r_2}$ and $B=\frac{r_2-r_0}{r_2-r_1}$. 
\begin{lem}\label{Lema:A,B,h,H}
The function $h_r$ is analytic and unimodal on the interval $\mathscr{I}_r$ with extremum at $u=r_0$. Moreover, if $u_0,u_1$ are the (potentially infinite) extreme points of $\mathscr{I}_r$ then $\lim\limits_{u\to u_0^+}h_r(u)=\lim\limits_{u\to u_1^-}h_r(u)$. In addition, the function
\begin{equation}\label{eq:H_r}
H_r(x,y)=h_r(y)-h_r(x),
\end{equation}
is a first integral for $\Omega_r$ on the fundamental square $\mathscr{I}_r^2$ and both separatrices of the saddle $\mathbf{s}$ are contained in the level set $\{H_r(x,y)=0\}$. Furthermore, the following holds:
\begin{itemize}
    \item[(i)] There exists a unique analytic strictly decreasing involution $\sigma_r:\mathscr{I}_r\to\mathscr{I}_r$ such that $h_r\circ\sigma_r=h_r$ with a unique fixed point $r_0$.
    \item[(ii)] The separatrix $\Gamma_r$ is given by the graph of $\sigma_r$ for $u\in\mathscr{I}_r$;
    \item[(iii)] If $(x,y)\in\Gamma_r$ and $x<y$, then $x<r_0<y$.
\end{itemize}
\end{lem}
\begin{proof}
The only points where $h_r$ fails to be analytic are $u=r_1$ and $u=r_2$. The fact that $h_r$ is unimodal follows easily from the computation
\begin{equation}\label{aux:1}
h_r'(u)=\frac{A}{u-r_1}+\frac{B}{u-r_2}=\frac{u-r_0}{(u-r_1)(u-r_2)}. 
\end{equation}
Now, computing the limits
\[\lim\limits_{|u|\to+\infty}h_r(u)=\lim\limits_{|u|\to +\infty}\log|u|+A\log\left\vert 1-\frac{r_1}{u}\right\vert+B\log\left\vert 1-\frac{r_2}{u}\right\vert=+\infty,\]
\[\lim\limits_{u\to r_1}\operatorname{sign}(r_0-r_1)h_r(u)=-\infty,\quad \lim\limits_{u\to r_2}\operatorname{sign}(r_2-r_0)h_r(u)=-\infty,\]
proves the first assertion. Again, by equation \eqref{aux:1}, we have that $\Omega_r=-dH_r$ on the fundamental square $\mathscr{I}_r^2$ and thus $H_r$ is a first integral for $\Omega_r$. Since $H_r(r_0,r_0)=0$, the saddle $\mathbf{s}$ along with both its separatrices are entirely contained in the level set $\{H_r(x,y)=0\}$. 

Now, we prove items (i)--(iii). Since $h_r$ is unimodal on the interval $\mathscr{I}_r$ with matching limits at the extreme points, Lemma \ref{Lema:unimodal-involution} implies that the equation $h_r(x)=h_r(y)$ defines an involution $y=\sigma_r(x)$ which proves (i).
To prove (ii), it is enough to observe that the level set $\{H_r= 0\}$ is defined precisely by $h_r(x)=h_r(y)$. Thus, either $x=y$ or $y=\sigma_r(x)$. By definition, $\Gamma_r$ must be the graph of $\sigma_r$.
Item (iii) follows from the fact that $\sigma_r$ exchanges $\{u\in\mathscr{I}_r :u<r_0\}$ with $\{u\in\mathscr{I}_r :u>r_0\}$.
\end{proof}

\begin{lem}\label{Lema:2r0=r1+r2}
$\Gamma_r$ is contained in a line if and only if $2 r_0=r_1+r_2$. In this case, $\Gamma_r$ is the graph of $\sigma_r(u)=2r_0-u$, for $u\in(r_1,r_2)$.
\end{lem}
\begin{proof}
We first remark that, if $\Gamma_r$ is contained in a line, then it must be the line $y=2r_0-x$ as we can see from \eqref{eq:DE3} that the non-diagonal eigenspace associated to $\mathbf{s}$ is spanned by the vector $(1,-1)$. If $2r_0=r_1+r_2$, by \eqref{eq:A,B,h}, for $u\in\mathscr{I}_r=(r_1,r_2)$, we have that $h_r(u)=\frac{1}{2}\log\left((u-r_1)(r_2-u)\right)$. The function $\tilde{\sigma}(u)=r_1+r_2-u$ is a non-trivial involution that satisfies
\[h_r\circ\tilde{\sigma}(u)=\frac{1}{2}\log\left((r_2-u)(u-r_1)\right)=h_r(u).\]
By Lemma \ref{Lema:A,B,h,H} (i), $\sigma_r=\tilde{\sigma}$. Conversely, if $\Gamma_r$ is given by $\sigma_r(u)=2r_0-u$, then $\tilde{H}(u):=H_r(u,2r_0-u)\equiv 0$. Since
\[\frac{d^3}{du^3}\tilde{H}(r_0)=\frac{4 (2 r_0-(r_1+r_2))}{(r_0-r_1)^2 (r_0-r_2)^2}=0,\]
we must have $2r_0=r_1+r_2$ and consequently $\mathscr{I}_r=(r_1,r_2)$.
\end{proof}

As we will see, our initial problem reduces to the study of the intersections between separatrices $\Gamma_1,\Gamma_2$ belonging to the family $\mathbf{\Gamma}$. More precisely, we consider the following notion:

\begin{defi}[Contact]
Given two smooth curves $\mathscr{C}_1,\mathscr{C}_2\subset\mathbb{R}^2$, we say that they have a \emph{contact of order $k$} (or a \emph{$k$-contact}) \emph{at the point $p=(a,b)$}, when the germs of $\mathscr{C}_1$ and $\mathscr{C}_2$ at $p$ coincide up to order $k-1$. Moreover, we say that the order of contact is \emph{exact} when $k$ is the greatest integer with such property.
\end{defi}

%\begin{figure}[!ht]
%    \centering
%    \begin{overpic}[width=0.25\linewidth]{2contact}
%    \put(49,49){$p$}
%    \end{overpic}
%    \caption{A $2$-contact between two curves.}
%    \label{Fig:contactdefi}
%\end{figure}

%\begin{defi}[Contact]
%Given two separatrices $\Gamma_1,\Gamma_2\in\mathbf{\Gamma}$, we say that they have a \emph{contact of order $k$} (or a \emph{$k$-contact}) \emph{at point $p=(a,b)$}, when the germs of $\Gamma_1$ and $\Gamma_2$ at $p$ coincide up to order $k$.
%\end{defi}

More concretely, since we can always write $\Gamma_1$ and $\Gamma_2$ as graphs of real analytic maps $y=\sigma_1(x)$ and $y=
\sigma_2(x)$, the above contact condition amounts to requiring that $\sigma_1(a)=\sigma_2(a)=b$ and
\[j^{k-1}(\sigma_1-\sigma_2)(a)\equiv 0,\]
where $j^mf(a)$ denotes the $m$-jet of a function $f(x)$ at $x=a$. If, additionally, we require the contact to be exact, we must have that $j^{k}(\sigma_1-\sigma_2)(a)\neq 0$. From the symmetry of the foliations, it follows immediately that contact points off the diagonal $x=y$ come in pairs, namely if $p=(a,b)$ is a $k$-contact point of $\Gamma_1$ and $\Gamma_2$, then its symmetric $q =(b, a)$ is also a $k$-contact point.

\subsection{First integrals for linear vector fields}\label{SS:FI}

An important fact which will be essential to the arguments in the present paper is that any diagonalizable linear vector field admits a Darboux first integral.  More precisely, the orbital Jordan form \eqref{eq:Jordan} admits $I(u,v)=\lambda\log|u|+\log|v|$ as a first integral. Thus, if $(u,v)=(\ell_1(x,y),\ell_2(x,y))$ is the affine change of variables that takes $X$ into its Jordan form, one can readily see that the function\begin{equation}\label{eq:DarbouxIntegral}
I(x,y)=\lambda\log|\ell_1(x,y)|+\log|\ell_2(x,y)|,
\end{equation}
is a first integral for $X$, well-defined in each one of the four quadrants of $\mathbb{R}^2\setminus\{\ell_1\ell_2=0\}$. Up to the absolute values, which are necessary since $\lambda$ need not be an integer, the function $I$ is the logarithm of the Darboux type function $|\ell_1|^\lambda|\ell_2|$. In the present work, it is more convenient to work with the logarithmic first integral. As we shall see in the next result, in the context studied in Proposition \ref{Propo:half-return}, we can fully characterize the half-return map of $X$ by considering the restriction of its first integral $I$ to the affine line $\Sigma$.

\begin{propo}\label{Propo:N(p)=N(q)}
Consider a half-plane $\mathcal{V}=\{f\geqslant0\}$ with boundary $\Sigma=\{f=0\}$ and a diagonalizable linear vector field $X$ such that $X\vert_\mathcal{V}$ has a well-defined half-return map $\sigma:\Sigma_{\rm in}\to\Sigma_{\rm out}$. Then the graph of $\sigma$ is defined by the relation
\[\{(p,q)\in\Sigma_{\rm in}\times\Sigma_{\rm out}:h(p)=h(q)\},\]
where $h=I\vert_\Sigma$ denotes the restriction of the first integral \eqref{eq:DarbouxIntegral} to $\Sigma$.
\end{propo}
\begin{proof}
In restriction to each quadrant, the level sets $\{I=c\}$ are smooth connected curves. Therefore each orbit of $X$ on a quadrant is in one-to-one correspondence with a level curve of $I$. As we have seen in Proposition \ref{Propo:half-return}, $\Sigma_{\rm in}$ and $\Sigma_{\rm out}$ are disjoint sets wholly contained in the same quadrant of $\mathbb{R}^2\setminus\{\ell_1\ell_2=0\}$. Moreover, each orbit starting at a point $p\in\Sigma_{\rm in}$ has a unique intersection with $\Sigma_{\rm out}$, which is precisely the point $q=\sigma(p)$. We conclude that $h(p)=h(q)$.

Conversely, if $p\in\Sigma_{\rm in}$ and $q\in\Sigma_{\rm out}$ are such that $h(p)=h(q)$, then $p$ and $q$ lie in the same orbit of $X$. Since any orbit of $X\vert_{\mathcal{V}}$ meets $\Sigma$ at most at two points, those must be exactly $p$ and $\sigma(p)$. The fact that $\Sigma_{\rm in}$ and $\Sigma_{\rm out}$ are disjoint implies that $q=\sigma(p)$.
\end{proof}

\begin{figure}[!ht]
    \centering
    \vspace*{0.5cm}
    \begin{overpic}[width=\linewidth]{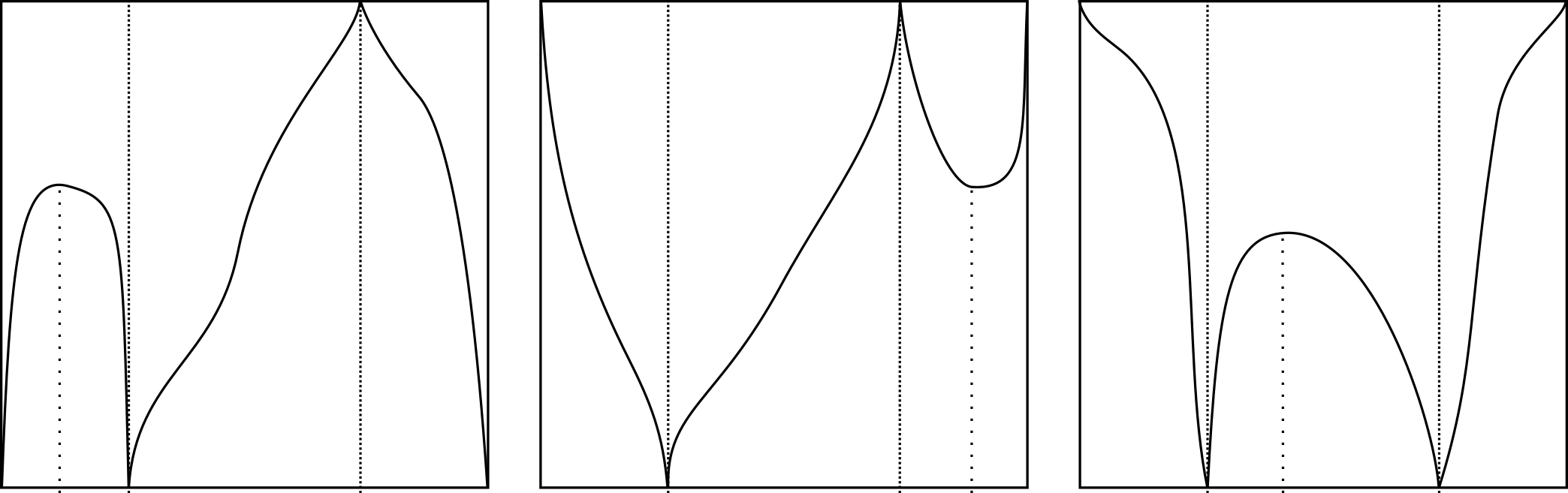}
    \put(9,32.5){(a) $\lambda<-1$}
    \put(3.5,-1.5){$r_0$}
    \put(8,-1.5){$r_1$}
    \put(23,-1.5){$r_2$}
    \put(40.5,32.5){(b) $-1<\lambda<0$}
    \put(62,-1.5){$r_0$}
    \put(42,-1.5){$r_1$}
    \put(57,-1.5){$r_2$}
    \put(79,32.5){(c) $\lambda>0$}
    \put(81.5,-1.5){$r_0$}
    \put(77,-1.5){$r_1$}
    \put(92,-1.5){$r_2$}
    \end{overpic}
    \caption{The graph of $h(x)$ for different values of $\lambda$ plotted on $\mathbb{RP}^1\times\mathbb{RP}^1$.}
    \label{Fig:Nunimodal}
\end{figure}

From this point forward, we establish the convention that $\Sigma=\{y=0\}$. Proposition \ref{Propo:half-return} implies that a necessary condition for $X$ to admit a well-defined half-return map on $\Sigma$ is that the invariant lines $\ell_1$ and $\ell_2$ are both transverse to $\Sigma$ at different points. Hence, up to additive constants, we can write the first integral as
\[I(x,y)=\lambda\log|x+\psi_1y-r_2|+\log|x+\psi_2y-r_1|,\]
where $\psi_1,\psi_2$ are distinct nonzero real numbers and $r_2$, $r_1$ are the abscissas of the intersections of $\Sigma$ with the lines $\ell_1$ and $\ell_2$, respectively. Since $I/\lambda$ is also a first integral, we may exchange the roles of $\ell_1$ and $\ell_2$, replacing $\lambda$ by $1/\lambda$, and assume the convention $r_1<r_2$. This is why both ranges $\lambda<-1$ and $-1<\lambda<0$ appear below. By Proposition \ref{Propo:N(p)=N(q)}, the half-return map can be defined by the level sets of the function
\begin{equation}\label{eq:N}
h(x)=I(x,0)=\lambda\log|x-r_2|+\log|x-r_1|,
\end{equation}
which has a unique critical point $r_0$, given by
\begin{equation}\label{eq:r0}
r_0=\frac{r_2+\lambda r_1}{1+\lambda}.
\end{equation}
Note that we have three distinct scenarios depending on $\lambda$:
\begin{itemize}
    \item[i.] If $\lambda<-1$, then $r_0\in\mathscr{I}= (-\infty,r_1)$;
    \item[ii.] If $-1<\lambda<0$, then $r_0\in\mathscr{I}=(r_2,+\infty)$;
    \item[iii.] If $\lambda>0$, then $r_0\in\mathscr{I}=(r_1,r_2)$. 
\end{itemize}
Figure~\ref{Fig:Nunimodal} displays the graph of $h(x)$ in each of these three scenarios.

Dynamically, the point $(r_0,0)$ is the invisible fold of $X$. Moreover, one can readily see by comparing equations \eqref{eq:N} and \eqref{eq:A,B,h} that
\[h(x)=(1+\lambda)h_r(x),\quad A=\frac{1}{1+\lambda}, B=\frac{\lambda}{1+\lambda},\]
for $r=(r_0,r_1,r_2)$. In particular, we have that $\mathscr{I}=\mathscr{I}_r$ as in \eqref{eq:Ir}. Since the factor $1+\lambda$ is nonzero, Lemma \ref{Lema:A,B,h,H} (i) implies that there exists a unique involution $\sigma:\mathscr{I}\to\mathscr{I}$ satisfying $h\circ\sigma(x)=h(x)$. Moreover, this involution is precisely $\sigma=\sigma_r$. Proposition \ref{Propo:N(p)=N(q)} implies that the half-return map is given by $(x,0)\mapsto(\sigma(x),0)$, for $x$ in a one-sided subinterval of $\mathscr{I}$ with $r_0$ as an extreme point. More precisely, we have that the half-return map is $(x,0)\mapsto(\sigma(x),0)$, with $\sigma:\Sigma_{\rm in}\to\Sigma_{\rm out}$, where, up to reorientation of the orbits,
\begin{itemize}
    \item[i.] $\Sigma_{\rm in}=(-\infty,r_0)$ and $\Sigma_{\rm out}=(r_0,r_1)$, for $\lambda<-1$;
    \item[ii.] $\Sigma_{\rm in}=(r_2,r_0)$ and $\Sigma_{\rm out}=(r_0,+\infty)$, for $-1<\lambda<0$;
    \item[iii.] $\Sigma_{\rm in}=(r_1,r_0)$ and $\Sigma_{\rm out}=(r_0,r_2)$, for $\lambda>0$.
\end{itemize}
Henceforth, we will denote the half-return map simply by its abscissa $x\mapsto\sigma(x)$. We remark that the intervals $\Sigma_{\rm in}$ and $\Sigma_{\rm out}$ above may be reversed according to the orientation of the orbits. Indeed, if $X$ has a well-defined half-return map $\sigma:\Sigma_{\rm in}\to\Sigma_{\rm out}$ constructed as above, then the same involution defines the half-return map $\sigma:\Sigma_{\rm out}\to\Sigma_{\rm in}$ of the vector field $-X$.

We are ready to state an essential fact about the half-return map.

\begin{propo}\label{Propo:Poincare-tripolar}
The graph of the involution $\sigma:\mathscr{I}_r\to\mathscr{I}_r$ is a solution curve of the foliation associated to the differential form $\Omega_r\in\mathbf{\Omega}$ in restriction to the square $\mathscr{I}^2$, passing through the point $\mathbf{s}=(r_0,r_0)$, with $r_0$ given by \eqref{eq:r0}. More precisely it coincides with the corresponding separatrix $\Gamma_r\in\mathbf{\Gamma}$. As a consequence the half-return map $\sigma:\Sigma_{\rm in}\to\Sigma_{\rm out}$ is one of the two branches of $\Gamma_r\setminus\{\mathbf{s}\}$.
\end{propo}
\begin{proof}
We have that 
\[\Omega_r=\frac{x-r_0}{(x-r_1)(x-r_2)}dx-\frac{y-r_0}{(y-r_1)(y-r_2)}dy=\frac{1}{1+\lambda}(h'(x)dx-h'(y)dy).\]
Taking $y=\sigma(x)$ for $x\in\mathscr{I}$, we have that $h(x)-h(y)\equiv 0$, and thus $h'(x)dx-h'(y)dy=0$. Hence, by the above equation $\Omega_r$ vanishes on $y=\sigma(x)$. The result follows from the fact that $\sigma(r_0)=r_0$ and from Lemma \ref{Lema:A,B,h,H} (ii), since $h(u)=(1+\lambda)h_r(u)$.
\end{proof}

\begin{figure}[!ht]
    \centering
    \begin{overpic}[width=0.6\linewidth]{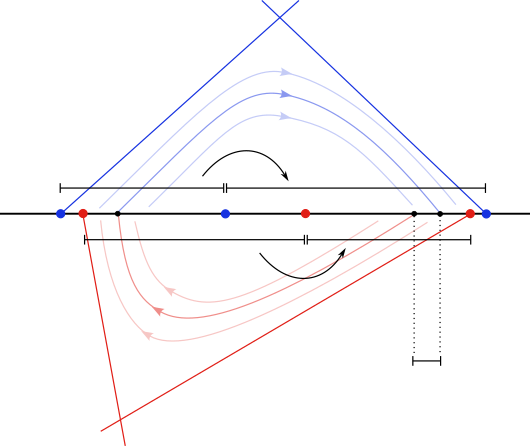}
    \put(20,45){$x$}
    \put(75,12.3){$d(x)=\sigma_+(x)-\sigma_{-}(x)$}
    \put(45.5,57){$\sigma_+$}
    \put(49.3,28){$\sigma_-^{-1}=\sigma_-$}
    \put(23,50.7){$\Sigma^+_{\rm in}$}
    \put(59,50.7){$\Sigma^+_{\rm out}$}
    \put(30,34.5){$\Sigma^{-}_{\rm out}$}
    \put(68.5,34.5){$\Sigma^{-}_{\rm in}$}
    \end{overpic}
    \caption{The half-return maps and the displacement function.}
    \label{Fig:displacement}
\end{figure}

Now, we consider a double diagonalizable piecewise linear vector field $X$ of the form \eqref{eq:1}, $\Sigma=\{y=0\}$ and $\Sigma^{\pm}=\{\pm y\geqslant0\}$. By Proposition \ref{Propo:half-return}, each component $X_\pm\vert_{\Sigma^{\pm}}$ has a well-defined half-return map $\sigma_\pm$ on $\Sigma$ if and only if it has an invisible fold $p^\pm$. We assume that this condition is satisfied and, additionally, that $X_{+}x(p^{+})\cdot X_{-}x(p^{-})<0$, so that the orbits of $X_+$ and $X_-$ turn in the same rotational sense. Then, we are able to define the displacement map $d:\Sigma_{\rm in}^+\cap\Sigma_{\rm out}^-\to\mathbb{R}$ by
\[d(x)=\sigma_{+}(x)-\sigma_{-}(x),\]
(see Figure \ref{Fig:displacement}).  Clearly, the limit cycles of $X$ are in one-to-one correspondence with the isolated zeros of $d$. Indeed, a limit cycle of $X$ intersects $\Sigma$ at two points with abscissas $x_0\in\Sigma^+_{\rm in}\cap\Sigma^-_{\rm out}$ and $x_1\in\Sigma^+_{\rm out}\cap\Sigma^-_{\rm in}$ such that necessarily $x_1=\sigma_+(x_0)$ and $\sigma_-(x_1)=x_0$. Thus, we can see that, on account of $\sigma_{\pm}$ being an involution
\[\sigma_-\circ\sigma_+(x_0)=x_0\iff \sigma_+(x_0)=\sigma_-(x_0)\iff d(x_0)=0.\]
Conversely, any isolated zero $x_0$ of $d(x)$ corresponds to the abscissa of an orbit such that $\sigma_+(x_0)=\sigma_-(x_0)=x_1$ and thus the arcs of $X_{\pm}\vert_{\Sigma^\pm}$ connecting $(x_0,0)$ to $(x_1,0)$ concatenate into an isolated closed orbit of $X$, thus a limit cycle. Accordingly, we say that a limit cycle of $X$ has \emph{multiplicity $k$} if its corresponding zero $x_0$ satisfies
\[j^{k-1}d(x_0)\equiv 0,\quad \text{and } j^kd(x_0)\neq 0.\]
Note that $X$ cannot have limit cycles if $\Sigma_{\rm in}^+\cap\Sigma_{\rm out}^-=\emptyset$.

The next result is crucial to establish the correspondence between limit cycles of $X$ and intersections of separatrices in $\mathbf{\Gamma}$.

\begin{corol}\label{Corol:crucial}
Let $X$ be a piecewise linear vector field of D-D type of the form \eqref{eq:1}, with well-defined half-return maps $\sigma_{\pm}:\Sigma^\pm_{\rm in}\to\Sigma^\pm_{\rm out}$. Then, there exist uniquely defined parameters $r^+,r^-\in\Lambda$ such that $\sigma_{\pm}=\sigma_{r^\pm}$ and each limit cycle of $X$ of multiplicity $k$ corresponds to a contact of order exactly $k$ between the separatrices $\Gamma_{r^+}$ and $\Gamma_{r^-}$ of the associated symmetric 1-forms $\Omega_{r^+}$ and $\Omega_{r^-}$.
\end{corol}
\begin{proof}
By Proposition \ref{Propo:Poincare-tripolar}, each half-return map $\sigma_{\pm}:\Sigma^\pm_{\rm in}\to\Sigma^\pm_{\rm out}$ is a solution curve of the foliation of the respective differential form
\[\Omega_{r^\pm}=\frac{x-r_0^\pm}{(x-r_1^\pm)(x-r_2^\pm)}dx-\frac{y-r_0^\pm}{(y-r_1^\pm)(y-r_2^\pm)}dy,\]
with $r_1^\pm<r_2^\pm$ and $r_0^\pm$ satisfying \eqref{eq:r0}, their graphs being the separatrices $\Gamma_{r^+}$ and $\Gamma_{r^-}$, respectively. Since the parameters $r_i^\pm$, $i=0,1,2$ are determined from the component $X_\pm$, the triples $r^\pm=(r_0^\pm,r_1^\pm,r_2^\pm)$ are uniquely determined by $X$. Let $x_0$ be an isolated zero of the displacement map. Then, the point $p=(x_0,\sigma_+(x_0))=(x_0,\sigma_-(x_0))$ is an intersection between $\Gamma_{r^+}$ and $\Gamma_{r^-}$. It corresponds to a limit cycle of multiplicity $k$ if and only if $j^{k-1}d(x_0)\equiv 0$, and $j^kd(x_0)\neq 0,$ which are precisely the conditions that $\Gamma_{r^+}$ and $\Gamma_{r^-}$ have a contact of order exactly $k$ at $p$.
\end{proof}

Corollary \ref{Corol:crucial} translates the investigation of the limit cycles of $X$ into the investigation of the intersections of the two separatrices $\Gamma_{r^+}$ and $\Gamma_{r^-}$ of Subsection~\ref{SS:Separatrices}. As previously observed, the existence of an intersection between $\Gamma_{r^+}$ and $\Gamma_{r^-}$ at a point $(x_0,x_1)$ implies that $(x_1,x_0)$ is also a point of intersection. Both intersections correspond to the same limit cycle, since only one point in the pair has abscissa located in the domain $\Sigma_{\rm in}^+\cap\Sigma_{\rm out}^-$.

The following result is the core of the present paper and from which Theorem \ref{Teo:A} is an immediate consequence in the piecewise linear vector fields scenario.

\begin{main}\label{Teo:I}
Let $\Gamma_1,\Gamma_2\in\mathbf{\Gamma}$ be two distinct separatrices satisfying one of the following conditions:
\begin{itemize}
    \item[1.] $\Gamma_1$ and $\Gamma_2$ have a $5$-contact at a diagonal point;
    \item[2.] $\Gamma_1$ and $\Gamma_2$ have a $3$-contact at a non-diagonal point.
\end{itemize}
Then $\Gamma_1$ and $\Gamma_2$ are nested segments of a line with equation $x+y=c$, one containing the other.
\end{main}

The order of contacts described in Theorem \ref{Teo:I} is optimal, that is, it is possible to find $\Gamma_1,\Gamma_2\in\mathbf{\Gamma}$ with a $4$-contact at a diagonal point and a $2$-contact at a non-diagonal point. We will discuss such scenarios in Examples \ref{Ex:1} and \ref{Ex:2} in Section \ref{S:Examples}. We now present a proof of Theorem \ref{Teo:A}.  Theorem \ref{Teo:I} will be proven in the next section.
\medskip

\noindent\emph{Proof of Theorem \ref{Teo:A}}: Let $X$ be a piecewise linear vector field \eqref{eq:1} of D-D type and suppose that it has a limit cycle of multiplicity $k\geqslant 3$. Let $(x_0,0)$ be an intersection between the limit cycle and $\Sigma$. By Corollary \ref{Corol:crucial}, there exist parameters $r^\pm$ such that $\Gamma_{r^+}$ and $\Gamma_{r^-}$ have a contact of order 3 at $p=(x_0,\sigma_{r^+}(x_0))$. Two cases must be considered. If $\sigma_{r^+}(x_0)=x_0$, then $\sigma_{r^-}(x_0)=x_0$ and thus $(x_0,0)$ is a tangential singularity and not a limit cycle. If $\sigma_{r^+}(x_0)\neq x_0$, then $p$ is a non-diagonal $3$-contact between the separatrices $\Gamma_{r^\pm}$. Thus, by Theorem \ref{Teo:I},  either the separatrices $\Gamma_{r^+}$ and $\Gamma_{r^-}$ coincide or are nested segments of the same line. Both cases imply that $d(x)\equiv 0$ in a neighborhood of $x_0$ which means that $x_0$ lies in a period annulus and cannot be the intersection of a limit cycle with $\Sigma$. This proves the first assertion. The sharpness of the bound is given by the example in Subsection \ref{Ex:2}.\qed

\section{Contact between separatrices in triple polar families}\label{S:contact}

In this section we will study the contact between two separatrices within the family $\mathbf{\Gamma}$ and prove Theorem \ref{Teo:I}. Cases 1 and 2 in Theorem \ref{Teo:I} are treated separately, as the tools we need to investigate the latter are substantially different from the ones used for the former. 

\subsection{Contact of order 5 at a diagonal point}\label{S:5contact}

Since the only intersection between any separatrix $\Gamma_r\in\mathbf{\Gamma}$ and the diagonal $\{x=y\}$ is the saddle point, two separatrices satisfying condition 1 in Theorem \ref{Teo:I} must have the $5$-contact at $\mathbf{s}=(r_0,r_0)$. The strategy to prove Theorem \ref{Teo:I} in this case is direct: we compare the $4$-jets of the respective functions $\sigma_r$ at the intersection point $x=r_0$.

Using Lemma \ref{Lema:A,B,h,H} (i), we can write
\begin{equation}\label{eq:pi_r}
\sigma_r(r_0+t)=r_0-t+\gamma_2t^2+\gamma_3t^3+\gamma_4t^4+\mathcal{O}(t^5).
\end{equation}
Thus, by Lemma \ref{Lema:A,B,h,H} (ii), the $4$-jet of $\Gamma_r$ at $\mathbf{s}$ is completely determined by the $5$-tuple $$(r_0,-1,\gamma_2,\gamma_3,\gamma_4).$$
In order to simplify the computations of $\gamma_k$, $k=2,3,4$, let us consider the nonzero quantities
\[\rho_1=r_0-r_1,\quad \rho_2=r_2-r_0.\]

\begin{propo}[$4$-jet of the separatrix at the saddle point]\label{Propo:4jet_pi_r}
Given a separatrix $\Gamma_r\in\mathbf{\Gamma}$, the expressions of the coefficients $\gamma_k$, $k=2,3,4$ in the expansion \eqref{eq:pi_r} of $\sigma_r(r_0+t)$ are given by
\[\gamma_2=\frac{2(\rho_2-\rho_1)}{3\rho_1\rho_2},\;\gamma_3=-\gamma_2^2,\;\gamma_4=\frac{11}{10}\gamma_2^3+\frac{\gamma_2}{5\rho_1\rho_2}.\]
\end{propo}
\begin{proof}
We first provide an expansion of the function $h_r$ as defined in \eqref{eq:A,B,h} for $|t|<\min\{|\rho_1|,|\rho_2|\}$. Setting $u=r_0+t$, we have that
\[|u-r_1|=|\rho_1+t|=|\rho_1|\left\vert 1+\frac{t}{\rho_1}\right\vert,\;|u-r_2|=|t-\rho_2|=|\rho_2|\left\vert 1-\frac{t}{\rho_2}\right\vert.\]
Thus,
\begin{align*}
h_r(r_0+t)&=A\log |\rho_1|+A\log\left\vert 1+\frac{t}{\rho_1}\right\vert+B\log|\rho_2|+B\log\left\vert 1-\frac{t}{\rho_2}\right\vert=\\
&=h_r(r_0)+A\log\left(1+\frac{t}{\rho_1}\right)+B\log\left(1-\frac{t}{\rho_2}\right)=\\
&=h_r(r_0)+\sum_{k\geqslant1}\left(\frac{(-1)^{k+1}A}{k\,\rho_1^k}-\frac{B}{k\,\rho_2^k}\right)t^k,
\end{align*}
where in the third equality, we use the power series expansion of $\log(1+t)$. Substituting $A=\rho_1/(\rho_1+\rho_2)$ and $B=\rho_2/(\rho_1+\rho_2)$, which hold by \eqref{eq:A,B,h}, yields
\begin{equation}\label{eq:h}
h_r(r_0+t)=h_r(r_0)+\sum_{k\geqslant2}h_k t^k,\quad h_k=-\frac{1}{k(\rho_1+\rho_2)}\left(\frac{1}{\rho_2^{k-1}}+\frac{(-1)^k}{\rho_1^{k-1}}\right).
\end{equation}
The absolute convergence of the power series of $\log(1+t)$ implies that the above series converges absolutely for $|t|<\min\{|\rho_1|,|\rho_2|\}$.
Now, let $\Delta(t):=\sigma_r(r_0+t)-r_0$. Since $h_r\circ\sigma_r=h_r$, we have the following equation
\[\sum_{k\geqslant 2}h_k\,\Delta(t)^k-\sum_{k\geqslant 2}h_k\;t^k=0.\]
Using the expansion \eqref{eq:pi_r} and equating the coefficients of $t^k$ in the above series to zero, we obtain the following set of equations:
\[2h_3+2h_2\gamma_2=0,\quad 3 h_3\gamma_2+h_2\gamma_2^2-2h_2\gamma_3=0,\;\]
\[ 2 h_2 \gamma_2 \gamma_3-2 h_5 - 4 h_4 \gamma_2 - 3 h_3 \gamma_2^2 + 3 h_3 \gamma_3  - 2 h_2 \gamma_4=0.\]
The result follows from substituting the expressions of $h_k$ given in equation \eqref{eq:h}.
\end{proof}

\begin{propo}\label{Propo:4jetunique}
Given $r\in\Lambda$, consider its corresponding 1-form $\Omega_r\in\mathbf{\Omega}$ and separatrix $\Gamma_r\in\mathbf{\Gamma}$ with expansion \eqref{eq:pi_r}. Then, the following holds:
\begin{itemize}
    \item[(i)] $\gamma_2=0$ if and only if $2r_0=r_1+r_2$. In that case, $\sigma_r(u)=r_1+r_2-u$, for $u\in(r_1,r_2)$;
    \item[(ii)] If $\gamma_2\neq 0$, then $r$ is uniquely determined by $(r_0,\gamma_2,\gamma_4)$, i.e., for every $\tilde{r}=(r_0,\tilde{r}_1,\tilde{r}_2)\in\Lambda$ such that $j^4(\sigma_{\tilde{r}}-\sigma_r)(r_0)=0$, then $\tilde{r}=r$.
\end{itemize}
\end{propo}
\begin{proof}
By Proposition \ref{Propo:4jet_pi_r}, $\gamma_2=0$ if and only if $\rho_1-\rho_2=2r_0-r_1-r_2=0$. In that case, Lemma \ref{Lema:2r0=r1+r2} yields assertion (i). Now, if $\gamma_2\neq 0$, then by Proposition \ref{Propo:4jet_pi_r}, we have that $10\gamma_4-11\gamma_2^3\neq 0$ and
\[\rho_1\rho_2=\frac{2\gamma_2}{10\gamma_4-11\gamma_2^3},\quad \rho_2-\rho_1=\frac{3\gamma_2^2}{10\gamma_4-11\gamma_2^3}.\]
Denoting by $\mathcal{P}$ and $\mathcal{D}$ the right-hand sides of the above equations, one can readily see that $-\rho_1$ and $\rho_2$ are the roots of the polynomial $x^2-\mathcal{D}x-\mathcal{P}$. Since the discriminant $\mathcal{D}^2+4\mathcal{P}=(\rho_1+\rho_2)^2$ is positive, the roots must be distinct and the following formulas hold:
\begin{equation}\label{aux:2}
r_1-r_0=-\rho_1=\frac{\mathcal{D}-\sqrt{\mathcal{D}^2+4\mathcal{P}}}{2},\quad r_2-r_0=\rho_2=\frac{\mathcal{D}+\sqrt{\mathcal{D}^2+4\mathcal{P}}}{2}.
\end{equation}
Now let $\tilde{r}=(r_0,\tilde{r}_1,\tilde{r}_2)$ and consider the corresponding involution \eqref{eq:pi_r}, denoted by 
\[\sigma_{\tilde{r}}(r_0+t)=r_0-t+\tilde{\gamma}_2t^2+\tilde{\gamma}_3t^3+\tilde{\gamma}_4t^4+\mathcal{O}(t^5).\] If $j^4(\sigma_{\tilde{r}}-\sigma_r)(r_0)=0$, then $\tilde{\gamma}_k=\gamma_k$ for $k=2,3,4$. Therefore, it follows by equation \eqref{aux:2}, that $\tilde{r}_1=r_1$ and $\tilde{r}_2=r_2$.
\end{proof}

\begin{obs}
Proposition \ref{Propo:4jetunique} states that the map $r\mapsto(r_0,\gamma_2(r),\gamma_4(r))$ is injective on $\{r\in\Lambda :2r_0\neq r_1+r_2\}$.  On the set $\{r\in\Lambda :2r_0 =  r_1+r_2\}$, the map takes on the value $(r_0,0,0)$ and the corresponding separatrix $\Gamma_r$ is given by the anti-diagonal straight line $x+y=2r_0$.
\end{obs}

\begin{teo}\label{Teo:5contact}
Let $\Gamma_{r^+},\Gamma_{r^-}\in\mathbf{\Gamma}$ be two distinct separatrices passing through a common diagonal point $p=(a,a)$. Then, either
\begin{itemize}
    \item[(i)] $p$ is an exact $2$-contact of $\Gamma_{r^+}$ and $\Gamma_{r^-}$;
    \item[(ii)] $p$ is an exact $4$-contact of $\Gamma_{r^+}$ and $\Gamma_{r^-}$;
    \item[(iii)] $\Gamma_{r^+}$ and $\Gamma_{r^-}$ are nested segments of the line $\{x+y=2 a\}$.
\end{itemize}
\end{teo}
\begin{proof}
Denoting $r^\pm=(r_0^\pm,r_1^\pm,r_2^\pm)$, since the intersection of any separatrix in $\mathbf{\Gamma}$ and the diagonal $\{x=y\}$ is the saddle point $\mathbf{s}$, we must have $r_0^-=r_0^+=a$. Keeping the superscript notation, by \eqref{eq:pi_r}, we can write
\[(\sigma_{r^+}-\sigma_{r^-})(a+t)=(\gamma_2^+-\gamma_2^-)t^2+(\gamma_3^+-\gamma_3^-)t^3+(\gamma_4^+-\gamma_4^-)t^4+\mathcal{O}(t^5).\]
Thus, the least exact order of contact of $\Gamma_{r^+}$ and $\Gamma_{r^-}$ at $p$ is 2. If $p$ is a $3$-contact, then $\gamma_2^+=\gamma_2^-$, which implies by Proposition \ref{Propo:4jet_pi_r}, $\gamma_3^+=\gamma_3^-$ and thus $p$ is in fact a $4$-contact. Finally, if $p$ is a $5$-contact, then $\gamma_k^+=\gamma_k^-$, for $k=2,3,4$. 

If $\gamma_2^+\neq 0$, by Proposition \ref{Propo:4jetunique} (ii), $r^+=r^-$ and $\Gamma_{r^+}$ and $\Gamma_{r^-}$ would coincide. Thus, $\gamma_2^+=0$. Then, Proposition \ref{Propo:4jetunique} (i) implies that $2r_0^\pm=r_1^\pm+r_2^\pm$, and $\sigma_{r^\pm}(u)=2a-u$, for $u\in(r_1^\pm,r_2^\pm)$. Hence $\Gamma_{r^\pm}$ are segments of the line $x+y=2a$, nested according to whether $r_1^+\leqslant r_1^-$ or $r_1^+\geqslant r_1^-$. 
\end{proof}

\begin{obs}
Theorem \ref{Teo:5contact} implies that a monodromic tangential singularity at the point $(r_0,0)$ is either a weak focus of order 1 (condition (i)), a weak focus of order 2 (condition (ii)) or a center (condition (iii) or coinciding separatrices). In the case of facing saddles, Art{\'e}s et al.~\cite{ArtesMAT} proved that at most two limit cycles bifurcate from such a point, the second one appearing when the perturbation creates a sliding segment.
\end{obs}

\subsection{Contact of order 3 at a non-diagonal point}\label{S:3contact}

We will now consider the intersection between two separatrices at a non-diagonal point $p=(a,b)$. The condition that the order of contact at $p$ is three is translated into the fact that they share a $2$-jet at the intersection. The strategy to prove Theorem \ref{Teo:I} in this case is the following: for a given $2$-jet, determine the subset of differential forms in $\mathbf{\Omega}$ such that their corresponding separatrices realize the $2$-jet at $p$. We shall prove that this subset either contains only one element or the separatrices are nested segments of the same straight line.\medskip

We recall that the phase space of the family $\mathbf{\Omega}$ is $\mathbb{M}=\mathbb{RP}^1\times\mathbb{RP}^1$, which is naturally equipped with the trivial projections on each of its two components. Let $J_x^2(\mathbb{M})$ denote the space of $2$-jets of curves in $\mathbb{M}$ which have a regular projection onto the first component. In other words, fixing a coordinate system $(x,y)$ in $\mathbb{M}$, an element of $J_x^2(\mathbb{M})$ can be identified with a $4$-tuple $J=(a,b,\nu,c)\in\mathbb{R}^4$ determined from a curve $y=f(x)$ by writing
\begin{equation}\label{eq:fjet}
f(a)=b,\quad f'(a)=\nu,\quad f''(a)=c.
\end{equation}
Given a $2$-jet $J\in J_x^2(\mathbb{M})$, let $\mathbf{\Omega}_J\subset\mathbf{\Omega}$ denote the subset of differential forms in the triple polar family that have a solution curve $y=f(x)$ satisfying \eqref{eq:fjet}. Denoting $R(u):=\frac{u-r_0}{(u-r_1)(u-r_2)}$, we can write the differential equation \eqref{eq:DE3} associated to $\Omega_r$ as 
\[\frac{dy}{dx}=\frac{R(x)}{R(y)},\]
for $y\neq r_0$. Thus, $\Omega_r\in\mathbf{\Omega}_J$ if and only if the following conditions hold:
\begin{equation}\label{eq:JetConditions}
y(a)=b,\quad \frac{R(a)}{R(b)}=\nu,\quad \frac{R'(a)-R'(b)\nu^2}{R(b)}=c.
\end{equation}
For brevity, we will call these equations the \emph{jet conditions}. Additionally, we need to consider the subset $\mathbf{\Omega}_{J,\text{sep}}\subset\mathbf{\Omega}_J$ of $1$-forms $\Omega_r$ such that the point $p=(a,b)$ lies on the corresponding separatrix $\Gamma_r$. One can readily characterize this subset by means of the first integral $H_r(x,y)$ defined in \eqref{eq:H_r}:
\begin{equation}\label{eq:SepConditions}
(a,b)\in\mathscr{I}_r^2,\quad\text{and } \,H_r(a,b)=0,
\end{equation}
where $\mathscr{I}_r^2$ is the fundamental square as defined in \eqref{eq:Ir}. These will be called the \emph{separatrix conditions}. We notice that, while the jet conditions are algebraic, the separatrix condition $H_r(a,b)=0$ determines a transcendental hypersurface in the parameter space $\Lambda$. 

%The set $\mathbf{\Omega}_{J,{\rm sep}}$ is precisely that of the symmetric differential forms \eqref{eq:Omega_r} such that the corresponding separatrices have a contact of order $3$ at the point $p=(a,b)$. 
Since we are interested in the case in which $p$ is non-diagonal, it is enough to consider $a<b$ due to the symmetry $(x,y)\mapsto (y,x)$.  We have that two separatrices $\Gamma_{1}, \Gamma_{2}\in\mathbf{\Gamma}$ have a contact of order $3$ at a point $p=(a,b)$ with $a<b$ if and only if the corresponding $1$-forms $\Omega_{1},\Omega_{2}$ are in the class $\mathbf{\Omega}_{J,{\rm sep}}$ for some common jet $J=(a,b,\nu,c)$, with $(\nu,c)\in\mathbb{R}^2$. 

We first apply a projective change of coordinates in order to reduce the number of free parameters of the original jet and separatrix conditions \eqref{eq:JetConditions} and \eqref{eq:SepConditions}. In order to do so, we introduce the following functions, defined for all $t>0$:
\begin{align}\label{eq:L,g}
L(t)&=\dfrac{t\log t}{t-1},\;\text{extended analytically to }L(1)=1,\nonumber\\
g(t;v)&=\dfrac{L(t)}{t+v}=\dfrac{t\log t}{(t-1)(t+v)},\;\text{for } v>0.
\end{align}
We then define the main auxiliary functions
\begin{align}\label{eq:V,K}
V(t_1,t_2)&=\frac{t_1L(t_2)-t_2L(t_1)}{L(t_1)-L(t_2)},\nonumber\\
K(t_1,t_2)&=\frac{t_1t_2}{V(t_1,t_2)}+t_1+t_2-\frac{V(t_1,t_2)(t_1+t_2+V(t_1,t_2))}{t_1t_2}.
\end{align}
Note that since $L'(t)=(t-1-\log t)/(t-1)^2>0$ for $t\neq 1$ and $L'(1)=1/2$, $L$ is a strictly increasing function, so $V$ is well-defined on the domain 
\[\mathbf{D}=\{(t_1,t_2)\in\mathbb{R}^2 :t_1,t_2>0,t_1\neq t_2\}.\]
Moreover, $V>0$ on $\mathbf{D}$. Indeed, for $t_1<t_2$ we have $L(t_1)<L(t_2)$ and, since $L(1/t)=L(t)/t$, then $t_2L(t_1)-t_1L(t_2)=t_1t_2\big(L(1/t_1)-L(1/t_2)\big)>0$. Hence $K$ is also well-defined on $\mathbf{D}$, and both $V$ and $K$ are invariant under the symmetry $(t_1,t_2)\mapsto(t_2,t_1)$. We are now ready to state the following result.

\begin{propo}\label{Propo:Reduction}
Let $J=(a,b,\nu,c)\in J_x^2(\mathbb{M})$, with $a<b$, and suppose that $\mathbf{\Omega}_{J,\text{sep}}\neq\emptyset$. Then, $\nu<0$. For $\Omega_r\in\mathbf{\Omega}_{J,\text{sep}}$, consider the corresponding $r=(r_0,r_1,r_2)\in\Lambda$ and define the projective coordinate change
\[\varphi(u)=\frac{u-b}{u-a}.\]
Consider the new variables $t_i = \varphi(r_i)$, for $i=0,1,2$. Then, $t_0=\frac{t_1t_2}{\nu}$, $t_1,t_2>0$, and the following equations are satisfied:
\begin{equation}\label{eq:JetSepConditions_simplified}
V(t_1,t_2)=v,\quad K(t_1,t_2)=k,
\end{equation}
where $v=-\nu$ and $k=\frac{c(a-b)}{\nu}+1+\nu$.
\end{propo}
\begin{proof}
Given $\Omega_r\in\mathbf{\Omega}_{J,\text{sep}}$, since $(a,b)\in\Gamma_r$, Lemma \ref{Lema:A,B,h,H} (iii) implies that $a<r_0<b$. Since $a,b\in\mathscr{I}_r$, they lie on the same side of $r_i$ for $i=1,2$ and therefore $R(a)$ and $R(b)$ have opposite signs, which by \eqref{eq:JetConditions} yields $\nu<0$.

Now, consider the change of variables $(z,w)=\Phi(x,y)$ given by
\[z=\frac{x-a}{x-b}=\frac{1}{\varphi(x)},\quad w=\frac{y-b}{y-a}=\varphi(y).\]
Again, since $a,b\in\mathscr{I}_r$, we obtain that $t_i=\varphi(r_i)>0$ for $i=1,2$. Now, let $y=f(x)$ be the curve satisfying \eqref{eq:fjet} and $w=\tilde{f}(z)$ describe the same curve on the coordinates $(z,w)$. Then, \eqref{eq:fjet} is translated into:
\begin{equation}\label{eq:ftildejet}
\tilde{f}(0)=0,\quad \tilde{f}'(0)=-\nu,\quad \tilde{f}''(0)=-(c(a-b)+2\nu(\nu+1)).
\end{equation}
Thus, we have 
\[\Phi^\ast\Omega_r=-\frac{(t_1-1)(t_2-1)}{t_0-1}\left(\frac{t_0z-1}{(z-1)(t_1z-1)(t_2z -1)}\,dz+\frac{w-t_0}{(w-1)(w-t_1)(w-t_2)}\,dw\right),\]
whose associated differential equation has the form
\begin{equation*}%\label{eq:dzdw}
\dfrac{dw}{dz}=-S(z)U(w),
\end{equation*}
with
\[S(z)=\frac{t_0z-1}{(z-1)(t_1z-1)(t_2z-1)},\quad U(w)=\frac{(w-1)(w-t_1)(w-t_2)}{w-t_0}.\]
Moreover, since
\[\dfrac{d^2w}{dz^2}=-S'(z)U(w)-S(z)U'(w)\frac{dw}{dz}=-S(z)U(w)\left(\frac{S'(z)}{S(z)}+\frac{U'(w)}{U(w)}\dfrac{dw}{dz}\right),\]
substituting into the jet conditions \eqref{eq:ftildejet} yields
\begin{equation}\label{eq:JetConditions_2}
\nu=\frac{t_1t_2}{t_0},\;\quad \frac{c(a-b)}{\nu}+1+\nu=-t_0+t_1+t_2+\frac{\nu}{t_1}+\frac{\nu}{t_2}-\frac{\nu}{t_0}.
\end{equation}
We hold onto this equation and will return to it after we deal with the separatrix condition \eqref{eq:SepConditions}. Using the coordinates $t_0,t_1,t_2$, we can write
\begin{align*}
H_r(a,b)=h_r(b)-h_r(a)&=A\log\left\vert \frac{b-r_1}{a-r_1}\right\vert+ B\log\left\vert\frac{b-r_2}{a-r_2}\right\vert=\\
&=\frac{r_1-r_0}{r_1-r_2}\log\left\vert \frac{b-r_1}{a-r_1}\right\vert+ \frac{r_2-r_0}{r_2-r_1}\log\left\vert\frac{b-r_2}{a-r_2}\right\vert=\\
&=\frac{(t_1-t_0)(t_2-1)}{(t_0-1)(t_1-t_2)}\log t_1+ \frac{(t_2-t_0)(t_1-1)}{(t_0-1)(t_2-t_1)}\log t_2=\\
&=\frac{(t_1-1)(t_2-1)}{(t_1-t_2)(t_0-1)}\left(\frac{t_1-t_0}{t_1-1}\log t_1-\frac{t_2-t_0}{t_2-1}\log t_2\right).
\end{align*}
Hence, using the first equality in \eqref{eq:JetConditions_2}, $H_r(a,b)=0$ if and only if 
\[\frac{t_1-t_0}{t_1-1}\log t_1-\frac{t_2-t_0}{t_2-1}\log t_2=\frac{t_1(\nu-t_2)}{\nu(t_1-1)}\log t_1-\frac{t_2(\nu-t_1)}{\nu(t_2-1)}\log t_2=0,\]
which is equivalent to $g(t_1;-\nu)=g(t_2;-\nu)$. To finish the proof, consider $v,k$ as given in \eqref{eq:JetSepConditions_simplified}. It suffices to prove that $g(t_1;v)=g(t_2;v)$ is equivalent to $V(t_1,t_2)=v$ and, under this equality, equation \eqref{eq:JetConditions_2} implies that $K(t_1,t_2)=k$.
In fact,
\[g(t_1;v)=g(t_2;v)\iff L(t_1)(t_2+v)=L(t_2)(t_1+v)\iff  v(L(t_1)-L(t_2))=t_1L(t_2)-t_2L(t_1), \]
which is equivalent to $V(t_1,t_2)=v$ since $L$ is strictly increasing. Now, substituting $t_0=t_1t_2/\nu$ and $\nu=-V(t_1,t_2)$ into \eqref{eq:JetConditions_2} yields precisely
\[k=\frac{t_1t_2}{V(t_1,t_2)}+t_1+t_2-\frac{V(t_1,t_2)}{t_1}-\frac{V(t_1,t_2)}{t_2}-\frac{V(t_1,t_2)^2}{t_1t_2}=K(t_1,t_2).\]
\end{proof}

The first equation in \eqref{eq:JetConditions_2} allows us to eliminate $t_0=t_1t_2/\nu$, so that the jet and separatrix conditions are reduced to the two equations \eqref{eq:JetSepConditions_simplified} in the two parameters $t_1$ and $t_2$. Thus, Proposition \ref{Propo:Reduction} establishes an injection between the family $\mathbf{\Omega}_{J,\rm sep}$ and the fibers of the map $\Theta:\mathbf{D}\to\mathbb{R}^2$ given by
\begin{equation}\label{eq:Theta}
\Theta(t_1,t_2)=(V(t_1,t_2),K(t_1,t_2)).    
\end{equation}
Therefore, part 2 of Theorem \ref{Teo:I} is reduced to the study of the fibers $\Theta^{-1}(v,k)$ of $\Theta$. If every fiber over the relevant parameters is at most one point (up to the symmetry $(t_1,t_2)\mapsto (t_2,t_1)$), then a common $2$-jet forces the separatrices to be the same. The next theorem states that the uniqueness assertion is true for every $(v, k)$ with $v > 0$ except one: the solution set of $\Theta(t_1,t_2)=(1,0)$ is an entire hyperbola. We shall see that this exception corresponds to the case that the corresponding separatrix $\Gamma_r$ is a line segment. 

\begin{figure}[!ht]
    \centering
    \begin{overpic}[width=0.85\linewidth]{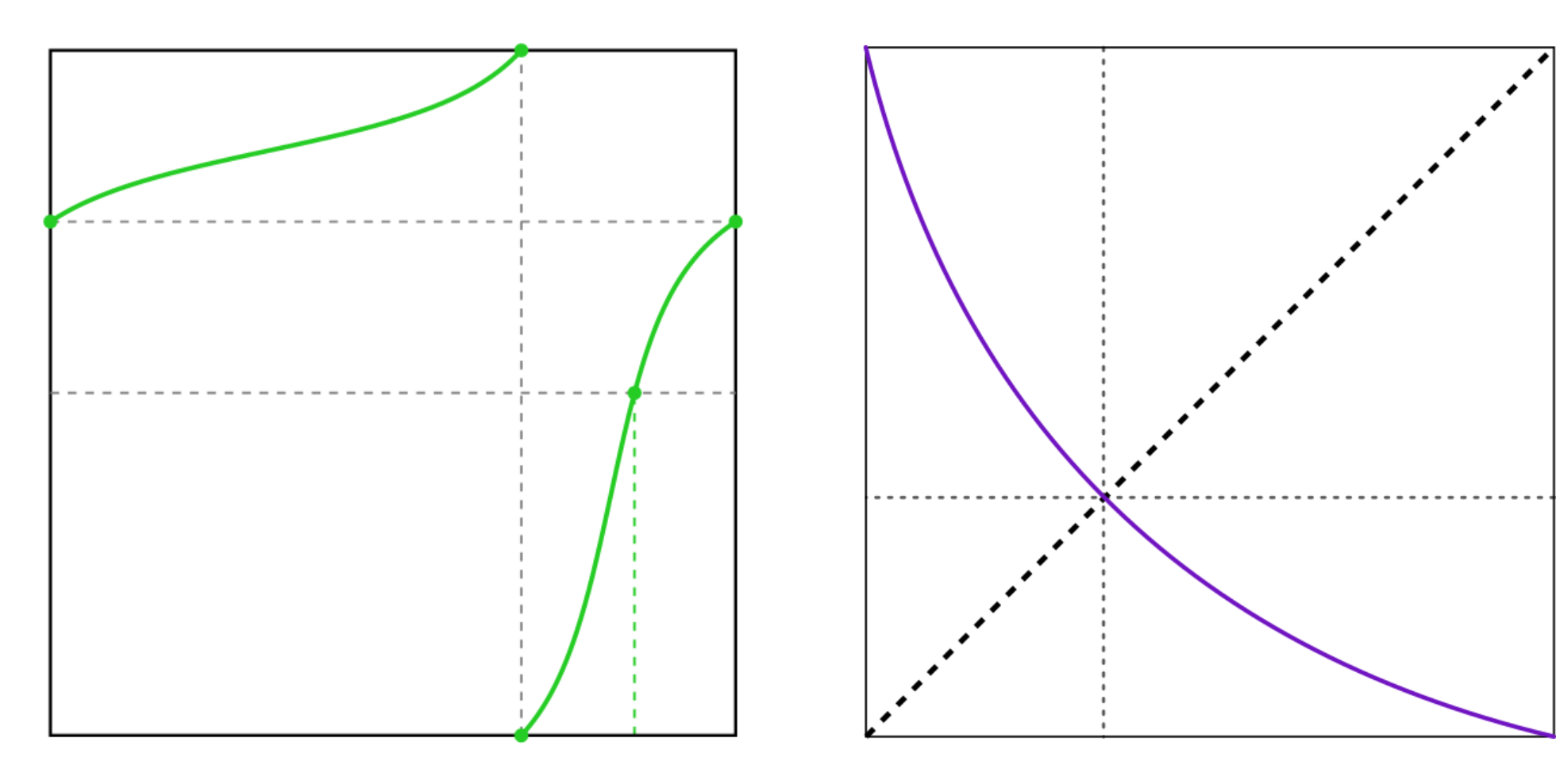}
    \put(32.7,0.5){$a$}
    \put(40,0.5){$b$}
    \put(48,36){$1$}
    \put(48,25){$0$}
    \put(48,46.5){$\infty$}
    \put(59,10){Nodes}
    \put(85,10){Saddles}
    \put(85.5,36){Nodes}
    \put(58.5,36){Saddles}
    \end{overpic}
    \caption{Left: The graph of $\varphi$ on $\mathbb{M}$. Right: A diagram representing each corresponding case in Remark \ref{Obs:translation} in the domain $\mathbf{D}$. The dashed thick black line represents the diagonal $t_1=t_2$ and the solid purple line represents the hyperbola $t_1t_2=1$. Colors available in the online version.}
    \label{Fig:varphi}
\end{figure}

\begin{obs}\label{Obs:translation}
The projective transformation $\varphi:\mathbb{RP}^1\to\mathbb{RP}^1$ sends the points $a,b,\infty$ to $\infty$, $0$ and $1$, respectively. Moreover, since $\varphi'(u)=\frac{b-a}{(u-a)^2}>0$, the map $\varphi$ is an increasing bijection, and $\varphi(u)>1$ for $u<a$, whereas $0<\varphi(u)<1$ for $u>b$. As $r_1,r_2\notin[a,b]$ and $r_0\in(a,b)$, we obtain
\begin{itemize}
\item[(i)] If $r_0<r_1$, then $0<t_1<t_2<1$;
\item[(ii)] If $r_0>r_2$, then $1<t_1<t_2$;
\item[(iii)] If $r_1<r_0<r_2$, then $0<t_2<1<t_1$.
\end{itemize}
In the D-D linear vector fields scenario, cases (i) and (ii) correspond to the nodes and case (iii) corresponds to the saddle (Figure \ref{Fig:varphi}).
\end{obs}

We consider the open domain
\[\mathbf{P}=\{(t_1,t_2)\in\mathbb{R}^2:0<t_1<t_2\},\]
and its decomposition $\mathbf{P}=\mathbf{Q}^+\sqcup\mathbf{H}\sqcup\mathbf{Q}^-$, where
\begin{equation*}%\label{eq:Q+HQ-}
\mathbf{Q}^+=\{(t_1,t_2)\in\mathbf{P} : t_1t_2>1\},\;\mathbf{H}=\{(t_1,t_2)\in\mathbf{P} : t_1t_2=1\},\;\mathbf{Q}^-=\{(t_1,t_2)\in\mathbf{P} : t_1t_2<1\}.
\end{equation*}
The map $\Theta$ is defined on the larger domain $\mathbf{D}$, yet the symmetry $(t_1,t_2)\mapsto (t_2,t_1)$ allows us to consider the subdomain $\mathbf{P}\subset\mathbf{D}$ which will make the computations easier. Now, we can state the following result. 

\begin{teo}\label{Teo:fibersTheta}
Consider the map $\Theta:\mathbf{D}\to\mathbb{R}^2$ given by \eqref{eq:Theta} and the decomposition $\mathbf{P}=\mathbf{Q}^+\sqcup\mathbf{H}\sqcup\mathbf{Q}^-$. Then, the following holds:
\begin{itemize}
\item[(a)] $\{(t_1,t_2)\in\mathbf{P}:V(t_1,t_2)=1\}=\mathbf{H}$ and $\Theta\vert_\mathbf{H}\equiv (1,0)$;
\item[(b)] $\Theta$ restricts to a diffeomorphism from $\mathbf{Q}^+$ onto its image, which is contained in $\{v>1\}$;
\item[(c)] $\Theta$ restricts to a diffeomorphism from $\mathbf{Q}^-$ onto its image, which is contained in $\{0<v<1\}$.
\end{itemize}
In particular, $\Theta^{-1}(1,0)=\{(t_1,t_2)\in\mathbf{D}:t_1t_2=1\}$, $\Theta^{-1}(1,k)=\emptyset$ for $k\neq 0$, and for every $(v,k)$ with $v\neq 1$ the fiber $\Theta^{-1}(v,k)$ is either empty or consists of exactly two points, namely $(t_1,t_2)$ and $(t_2,t_1)$ for a unique $(t_1,t_2)\in\mathbf P$.
\end{teo}

\begin{figure}[!ht]
    \centering
    \begin{overpic}[width=0.85\linewidth]{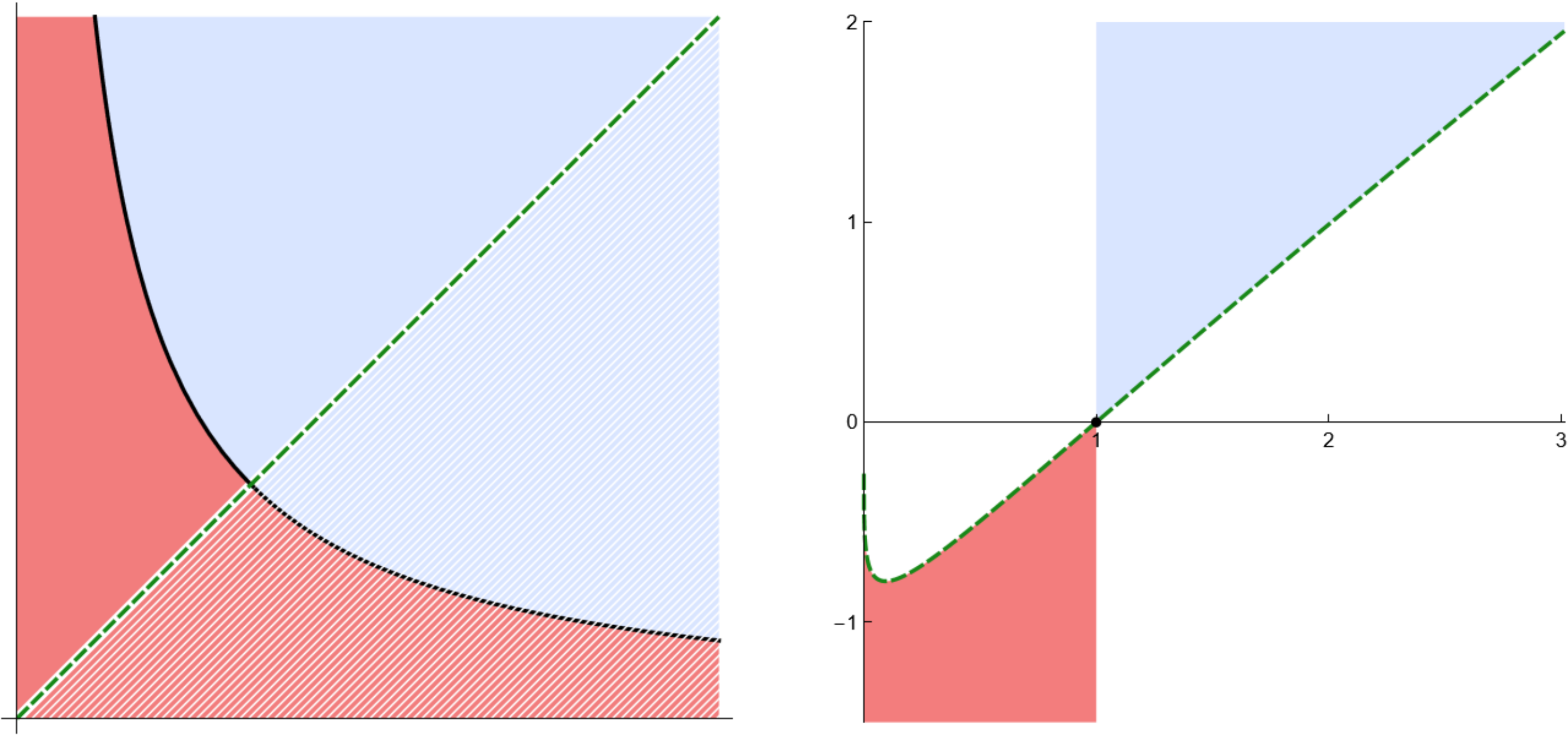}
    \put(5,15){$\mathbf{Q}^{-}$}
    \put(11,26){${\mathbf{H}}$}
    \put(23,35){$\mathbf{Q}^+$}
    \put(47,1){$t_1$}
    \put(1,47){$t_2$}
    \put(58,7){$\Theta(\mathbf{Q}^-)$}
    \put(75,35){$\Theta(\mathbf{Q}^+)$}
    \put(62,21){$\Theta({\mathbf{H}})$}
    \put(100.5,20){$v$}
    \put(55,46){$k$}
    \end{overpic}
    \caption{The mapping $\Theta:\mathbf{D}\to\mathbb{R}^2$. The subdomains $\mathbf{Q}^-$, $\mathbf{H}$ and $\mathbf{Q}^+$ and their reflections with respect to the diagonal $t_1=t_2$, represented in the hatched region with the same shading, have the same images, as $\Theta(t_1,t_2)=\Theta(t_2,t_1)$. The diagonal $t_1=t_2$ is displayed in dashed green line, together with the curve formed by the limits of $\Theta(t_1,t_2)$ as $(t_1,t_2)$ approaches the diagonal from either side. The images $\Theta(\mathbf{Q}^+)$ and $\Theta(\mathbf{Q}^-)$ lie on opposite sides of this curve with respect to $\Theta(\mathbf{H})=(1,0)$. Colors available in the online version.}
    \label{Fig:ImageTheta}
\end{figure}

The proof of Theorem \ref{Teo:fibersTheta} is very technical and relies on an extensive chain of results which would make a long diversion if presented here. For the sake of readability, we postpone the proof to Section \ref{S:Proofs} and finish the reasoning of the present section to establish the proof of Theorem \ref{Teo:I}. A graphical representation of the domain and image of the mapping $\Theta$ is given in Figure \ref{Fig:ImageTheta}. We carry on with the following results.

\begin{corol}\label{Corol:Omega_J,sep}
Given $J=(a,b,\nu,c)\in J_x^2(\mathbb{M})$, with $a<b$ and $\nu\neq -1$, the family $\mathbf{\Omega}_{J,\text{sep}}$ contains at most one element.
\end{corol}
\begin{proof}
By Proposition~\ref{Propo:Reduction}, the map $\Omega_r\mapsto(t_1,t_2)$ is an
injection of $\mathbf{\Omega}_{J,\mathrm{sep}}$ into $\Theta^{-1}(v,k)$. Indeed, $r=(r_0,r_1,r_2)$ is recovered from $(t_0,t_1,t_2)$ through $t_0=t_1t_2/\nu$ and $r_i=\varphi^{-1}(t_i)$. Since $v=-\nu\neq1$, Theorem \ref{Teo:fibersTheta} implies that the fiber $\Theta^{-1}(v,k)$ is either empty or contains exactly two points, namely $(t_1,t_2)$ and $(t_2,t_1)$ if $(v,k)\in\Theta(\mathbf Q^+)\cup\Theta(\mathbf Q^-)$. In the case that the fiber is not empty, suppose that both points are attained for two distinct parameter values $r$ and $\tilde{r}$. Thus $t_1=\tilde{t}_2$ and $t_2=\tilde{t}_1$. Then $\tilde{r}_1=\varphi^{-1}(\tilde{t}_1)=r_2>r_1=\varphi^{-1}(\tilde{t}_2)=\tilde{r}_2$ and thus $\tilde{r}\not\in\Lambda$. Hence $\mathbf\Omega_{J,\mathrm{sep}}$ has at most one element.    
\end{proof} 

\begin{teo}\label{Teo:3contact}
Given a jet $J=(a,b,\nu,c)$, with $a<b$, let $\Omega_{r^+},\Omega_{r^-}\in\mathbf{\Omega}_{J,{\rm sep}}$ and $\Gamma_{r^+}$ and $\Gamma_{r^-}$ be their respective separatrices. Then, the following holds:
\begin{itemize}
    \item[(i)] If $(\nu,c)\neq(-1,0)$, then $r^{+}=r^{-}$;
    \item[(ii)] If $(\nu,c)=(-1,0)$, then $\Gamma_{r^{+}}$ and $\Gamma_{r^{-}}$ are nested segments of the line $\{x+y=a+b\}$.
\end{itemize}
\end{teo}
\begin{proof}
Again, we retrieve the expressions for $v,k$ in \eqref{eq:JetSepConditions_simplified}, to conclude that $(\nu,c)=(-1,0)$ if and only if $(v,k)=(1,0)$. Item (i) follows immediately from Corollary \ref{Corol:Omega_J,sep} and Theorem \ref{Teo:fibersTheta}. Now, consider $(\nu,c)=(-1,0)$. For every $\Omega_r\in\mathbf{\Omega}_{J,{\rm sep}}$, Proposition \ref{Propo:Reduction} implies that the corresponding quantities $t_i$, $i=0,1,2$, satisfy $t_0=-t_1t_2$ and $(t_1,t_2)\in\Theta^{-1}(1,0)$. Thus, by Theorem \ref{Teo:fibersTheta}, $t_1t_2=1$ and we can write $(t_1,t_2)=(t,1/t)$ for some $t\neq 1$. Hence,
\[r_0=\varphi^{-1}(-1)=\frac{a+b}{2},\;r_1=\varphi^{-1}(t)=a-\frac{b-a}{t-1},\;r_2=\varphi^{-1}(1/t)=b+\frac{b-a}{t-1},\]
and thus $2 r_0=r_1+r_2$. Therefore, by Lemma \ref{Lema:2r0=r1+r2}, $\Gamma_r$ is contained in the line $\{x+y=a+b\}$. The injectivity of the function $(b-a)/(t-1)$ for $t\neq 1$ proves the assertion that any two such separatrices must be nested.
\end{proof}

Theorem \ref{Teo:I} is a direct consequence of Theorems \ref{Teo:5contact} and \ref{Teo:3contact}.\medskip

\noindent\emph{Proof of Theorem \ref{Teo:I}:} Let $\Gamma_1,\Gamma_2\in\mathbf{\Gamma}$ be distinct separatrices. Suppose that $\Gamma_1$ and $\Gamma_2$ have a 5-contact at a diagonal point $p=(a,a)$. Then, by Theorem \ref{Teo:5contact}, we have that $\Gamma_1$ and $\Gamma_2$ are nested segments of the line $\{x+y=2a\}$. Now, suppose that $\Gamma_1$ and $\Gamma_2$ have a 3-contact at a non-diagonal point $p=(a,b)$. Up to the symmetry $(x,y)\mapsto (y,x)$, we may assume that $a<b$. Thus, the respective 1-forms $\Omega_1,\Omega_2$ are such that $\Omega_1,\Omega_2\in\mathbf{\Omega}_{J,{\rm sep}}$ for some common jet $J=(a,b,\nu,c)$. By Theorem \ref{Teo:3contact}, we must have $(\nu,c)=(-1,0)$, otherwise $\Gamma_1=\Gamma_2$. Item (ii) of Theorem \ref{Teo:3contact} implies that $\Gamma_1$ and $\Gamma_2$ are nested segments of the line $\{x+y=a+b\}$. This completes the proof. \qed 

\section{Examples}\label{S:Examples}

\subsection{A 4-contact at a diagonal point}\label{Ex:1}
Let us fix the parameter value $r^+=(-1/3,-1,1)$ and its associated differential form and corresponding first integral defined on $\mathscr{I}_{r^+}^2=(-1,1)\times(-1,1)$:
\begin{equation}\label{eq:ex4contact1}
\Omega_{r^+}=\frac{3x+1}{3(x^2-1)}dx-\frac{3y+1}{3(y^2-1)}dy,\quad H_{r^+}(x,y)=\frac{1}{3} \log \left(\frac{(y-1)^2 (y+1)}{(x-1)^2 (x+1)}\right).
\end{equation}
The level curve $\{H_{r^+}=0\}$ is given by the polynomial equation $(x-y)F_1(x,y)=0$, with $F_1(x,y)=x^2+x y-x+y^2-y-1$. Thus, the separatrix $\Gamma_{r^+}$ is given by the intersection of the conic $\{F_1(x,y)=0\}$ and the square $(-1,1)\times(-1,1)$.
Similarly, for the parameter values $r^-=(-1/3,-11/9,7/3)$, we have
\begin{equation}\label{eq:ex4contact2}
\Omega_{r^-}=\frac{27x+9}{(9x+11)(3x-7)}dx-\frac{27y+9}{(9y+11)(3y-7)}dy,\quad H_{r^-}(x,y)=\frac{1}{4} \log \left(\frac{(3 y-7)^3 (9 y+11)}{(3 x-7)^3 (9 x+11)}\right),
\end{equation}
and we obtain that the separatrix $\Gamma_{r^-}$ is defined by the intersection of the square $(-11/9,7/3)\times(-11/9,7/3)$ with the algebraic curve $\{F_2(x,y)=0\}$, where 
$$F_2(x,y)=27 x^3+27 x^2 y-156 x^2+27 x y^2-156 x y+210 x+27 y^3-156 y^2+210 y+196.$$
Using the explicit rational parametrization of $\{F_1(x,y)=0\}$ given by
\[x(t)=\frac{t^2-t-1}{t^2+t+1},\quad y(t)=-\frac{t^2+t-1}{t^2+t+1},\]
which satisfies $(x(1),y(1))=\mathbf{s}=(-1/3,-1/3)$, we obtain
$F_2(x(t),y(t))=\frac{16}{3}(t-1)^4+O\big((t-1)^5\big)$. Since $\nabla F_2(\mathbf{s})=(384,384)\neq 0$, the separatrices $\Gamma_{r^+}$ and $\Gamma_{r^-}$ have an exact contact of order $4$ at $\mathbf{s}$.

\begin{figure}[!ht]
    \centering
    \begin{overpic}[width=0.85\linewidth]{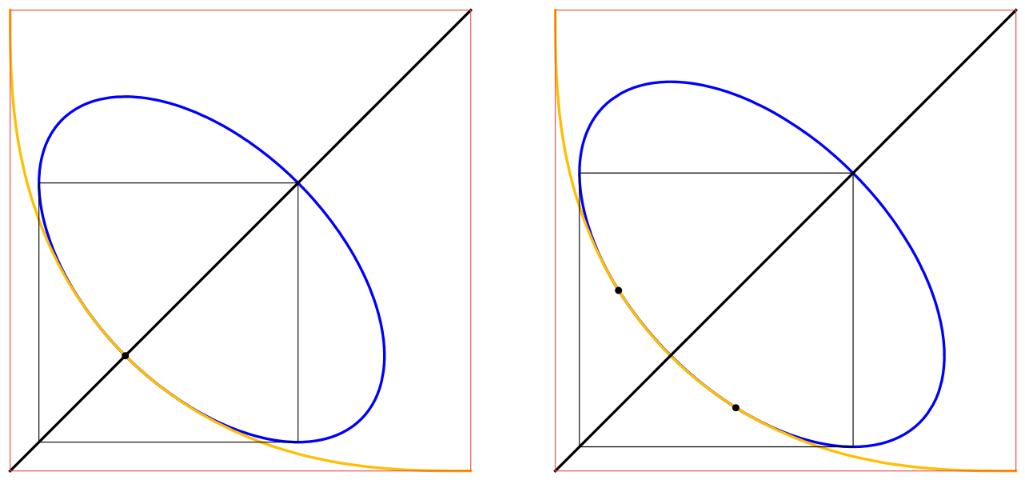}
    \put(11.8,13.5){$\mathbf{s}$}
    %\put(35,22){$\Gamma_{r^+}$}
    %\put(7,70){$\Gamma_{r^-}$}
    \put(58,16.5){$p$}
    %\put(35,13){$p$}
    %\put(35,22){$\Gamma_{r^+}$}
    %\put(5,85){$\Gamma_{r^\ast}$}
    \end{overpic}
    \caption{To the left: The $4$-contact between $\Gamma_{r^+}$ and $\Gamma_{r^-}$ for $r^+=(-1/3,-1,1)$ and $r^-=(-1/3,-11/9,7/3)$ at the common saddle $\mathbf{s}=(-1/3,-1/3)$. The conic $\{F_1(x,y)=0\}$ is drawn in blue and the relevant arc of the curve $\{F_2(x,y)=0\}$ is drawn in yellow. To the right: The $2$-contact between $\Gamma_{r^+}$ and $\Gamma_{r^\ast}$ at $p=(-5/7,1/7)$. The conic $\{F_1(x,y)=0\}$ is again drawn in blue and the relevant arc of the curve $\{F_2((x-\zeta)/\kappa,(y-\zeta)/\kappa)=0\}$ is drawn in yellow. Colors available in the online version.}
    \label{Fig:ex}
\end{figure}

\subsection{A 2-contact at a non-diagonal point}\label{Ex:2}
\begin{figure}[!ht]
    \centering
    \includegraphics[width=0.85\linewidth]{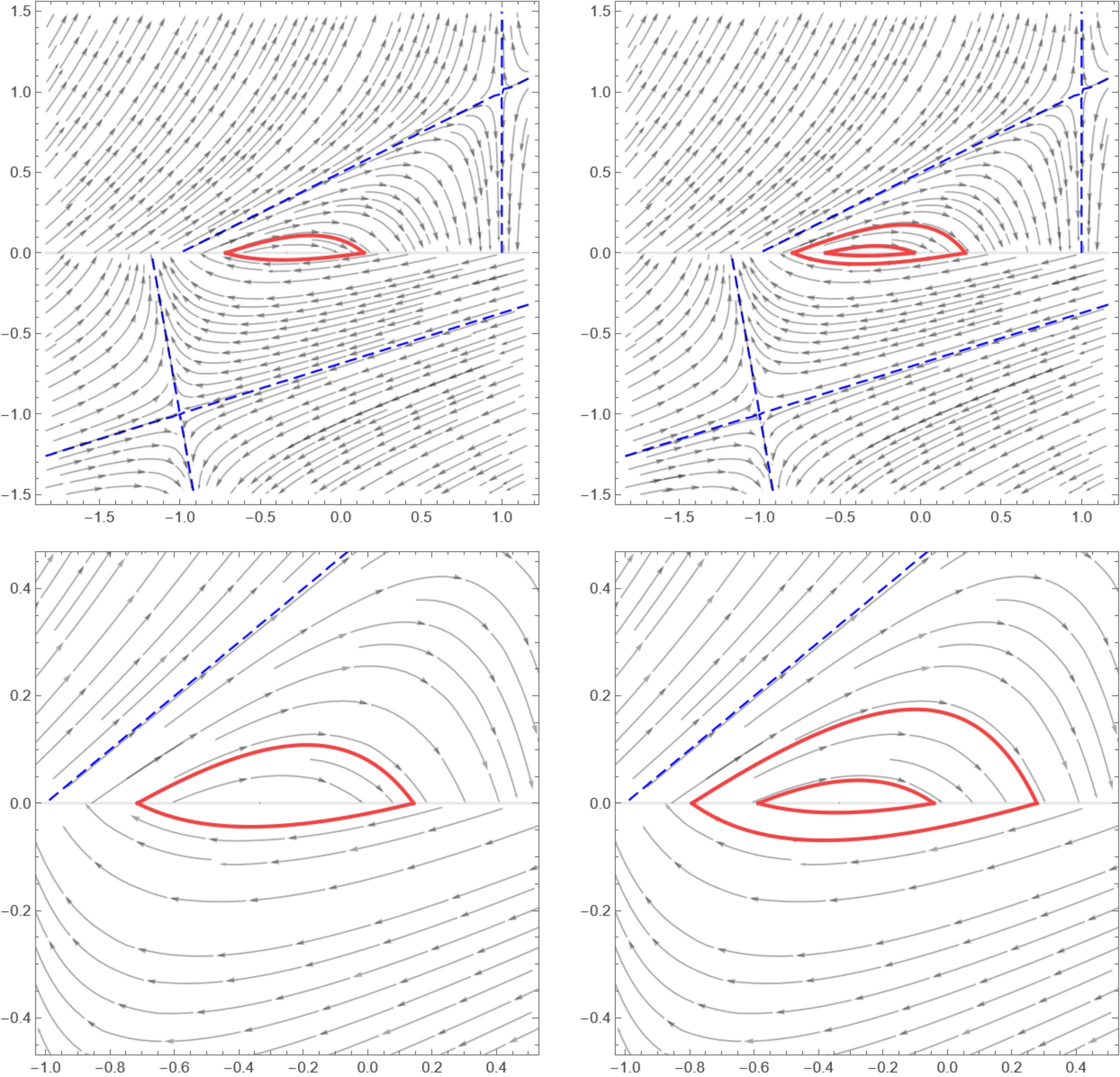}
    \caption{Left: The phase portrait of a piecewise linear vector field whose components correspond to  the parameters $r^{+}$ and $r^\ast$ in Example \ref{Ex:2}. The Poincaré half-maps correspond to parametrizations of the separatrices $\Gamma_{r^+}$ and $\Gamma_{r^\ast}$ and the $2$-contact at $p=(-5/7,1/7)$ corresponds to a limit cycle of multiplicity $2$ (in red) whose intersections with $\Sigma$ are $x=-5/7$ and $x=1/7$. Right: The limit cycle unfolds into two hyperbolic limit cycles under the perturbation $r(\varepsilon)=(r_0^\ast+\varepsilon,r_1^\ast+\varepsilon,r_2^\ast+\varepsilon)$ for $\varepsilon=5\times10^{-4}$.}
    \label{Fig:exDD}
\end{figure}

A slight modification on the previous example shows that distinct separatrices can have a $2$-contact at non-diagonal points. Consider $\Omega_{r^\pm}$ as in \eqref{eq:ex4contact1} and \eqref{eq:ex4contact2} along with the polynomials $F_1,F_2$ and the point $p=(-5/7,1/7)$, which satisfies $F_1(p)=0$, i.e. $p\in\Gamma_{r^+}$. For $\kappa>0$ and $\zeta\in\mathbb{R}$, let us consider the affine coordinate change $\phi(x,y)=((x-\zeta)/\kappa,(y-\zeta)/\kappa)$. We claim that there exist values of $(\kappa,\zeta)$ such that
\[F_2(\phi(p))=0,\;\text{and } \nabla F_1(p)\wedge\nabla F_2(\phi(p))=0.\]
As a consequence, the transformed differential form $\Omega_{r^\ast}=\phi^\ast\Omega_{r^-}$, with $r^\ast=(r_0^\ast,r_1^\ast,r_2^\ast)$, $r_i^\ast=\kappa r_i^-+\zeta$, has its corresponding separatrix $\Gamma_{r^\ast}$ with a contact of order 2 with $\Gamma_{r^+}$ at $p$. In order to prove that such values of $(\kappa,\zeta)$ exist, we write the expressions above as rational functions of $(\kappa,\zeta)$:
\[F_2(\phi(p))=\frac{4P_1(\kappa,\zeta)}{343\kappa^3},\;\nabla F_1(p)\wedge\nabla F_2(\phi(p))=\frac{36P_2(\kappa,\zeta)}{343\kappa^2}.\]
The ideal $\langle P_1,P_2\rangle$ admits a Gr{\"o}ebner basis $\{Q_1(\zeta),\kappa-Q_2(\zeta)\}$, where
\begin{align*}
Q_1(\zeta)&=413048832 \zeta ^4+220924928 \zeta ^3+24267936 \zeta ^2+5374656 \zeta +95461,\\
Q_2(\zeta)&=\frac{4988196864}{51075643}\zeta ^3+\frac{2559078144}{51075643}\zeta ^2+\frac{201982833  }{51075643}\zeta+\frac{359200557}{357529501}.
\end{align*}
Since, $Q_1(-1/50)Q_1(0)<0$, by continuity there exists $\zeta_0\in(-1/50,0)$ such that $Q_1(\zeta_0)=0$. Taking $(\kappa,\zeta)=(Q_2(\zeta_0),\zeta_0)$ proves the claim. A numerical approximation for such values is $\kappa\approx 0.946659442252134$ and $\zeta\approx -0.0191371680316586$, which correspond to the parameter values 
\[r^\ast=(r_0^\ast,r_1^\ast,r_2^\ast)\approx(-0.334690315449037, -1.17616537522871, 2.18973486388999).\]
The contact between $\Gamma_{r^+}$ and $\Gamma_{r^\ast}$ at $p$ is exactly $2$. Indeed, since $r_0^\ast\neq r_0^+$, the corresponding saddle points of the separatrices are distinct and thus they cannot coincide. By Theorem \ref{Teo:3contact}, if $p$ was a $3$-contact, then $\Gamma_{r^+}$ would be a segment of the line $x+y=-4/7$. Since $2 r_0^+\neq0= r_1^++r_2^+$, by Lemma \ref{Lema:2r0=r1+r2} this is not the case. An explicit example of a D-D vector field \eqref{eq:1} realizing this $2$-contact, that is, having a limit cycle with multiplicity $2$ is the following
\begin{align*}
X_+&=(1-x)\partial_x-\frac{1}{2}(1+3x-4y)\partial_y,\\
X_-&=(r_1^\ast-r_2^\ast)\biggl(\bigl((4+r_1^\ast+3r_2^\ast)x-4 (r_1^\ast+1) (r_2^\ast+1) y-4 r_1^\ast r_2^\ast-3 r_1^\ast-r_2^\ast\bigr)\partial_x+\\
& \qquad\qquad\qquad\qquad \bigl(4 x-3 r_1^\ast
   (y+1)-(r_2^\ast+4) y-r_2^\ast\bigr)\partial_y\biggr).
\end{align*}
Figures \ref{Fig:ex} and \ref{Fig:exDD} illustrate the example.

Under small perturbations, two $1$-contacts can unfold from the $2$-contact of $\Gamma_{r^+}$ and $\Gamma_{r^\ast}$. Consider, for instance, $r(\varepsilon)=(r_0^\ast+\varepsilon,r_1^\ast+\varepsilon,r_2^\ast+\varepsilon)$ and let $\sigma_{r^\ast}(u)$ parametrize the separatrix $\Gamma_{r^\ast}$, i.e. $h_{r^\ast}\circ\sigma_{r^\ast}=h_{r^\ast}$. Then, we have, by \eqref{eq:A,B,h},  $h_{r(\varepsilon)}(u)=h_{r^\ast}(u-\varepsilon)$ and the function $\sigma_{r(\varepsilon)}(u)=\sigma_{r^\ast}(u-\varepsilon)+\varepsilon$ satisfies $h_{r(\varepsilon)}\circ\sigma_{r(\varepsilon)}(u)=h_{r(\varepsilon)}(u)$, and thus $\sigma_{r(\varepsilon)}$ parametrizes $\Gamma_{r(\varepsilon)}$. Therefore, the displacement function is given by
\[d(x)=\sigma_{r(\varepsilon)}(x)-\sigma_{r^+}(x)=\sigma_{r^\ast}(x-\varepsilon)-\sigma_{r^+}(x)+\varepsilon.\]
Writing $x_0=-5/7$ and $d_0=\sigma_{r^\ast}-\sigma_{r^+}$, we have $d_0(x_0)=d_0'(x_0)=0$, $d_0''(x_0)=\sigma_{r^\ast}''(x_0)-\sigma_{r^+}''(x_0)\approx -0.286418<0$ and $\partial_\varepsilon d(x_0)\vert_{\varepsilon=0}=1-\sigma_{r^\ast}'(x_0)=13/5>0$. Hence, for small $\varepsilon>0$, $d$ has two simple zeros close to $x_0$ and none for $\varepsilon<0$. In other words, the $2$-contact at $\varepsilon=0$ becomes two contacts of order $1$. The unfolding is depicted in Figure \ref{Fig:exDD}.

\section{Technical lemmas and Proof of Theorem \ref{Teo:fibersTheta}}\label{S:Proofs}

This section will be dedicated to prove Theorem \ref{Teo:fibersTheta}. The proof rests on three facts that will be established in the following subsections. First, we prove that the level curves $\mathcal{C}_v=\{V(t_1,t_2)=v\}$ are path-connected and wholly contained in a single subdomain of the decomposition $\mathbf{P}=\mathbf{Q}^+\sqcup\mathbf{H}\sqcup\mathbf{Q}^-$. Thus, we reduce the study of the map $\Theta$ to that of the restriction $K\vert_{\mathcal{C}_v}$ of the function $K$ to the level curves $\mathcal{C}_v$. Second, we will prove that $\det D\Theta$ does not vanish on $\mathbf{Q}^{\pm}$. Hence, $K\vert_{\mathcal{C}_v}$ is monotone and therefore $\Theta$ must be injective and thus a diffeomorphism on each subdomain $\mathbf{Q}^\pm$. Finally, we compute the limit points of $K\vert_{\mathcal{C}_v}$ to determine precisely the image of $\Theta$ on each subdomain. 

\subsection{The level curves \texorpdfstring{$V(t_1,t_2)=v$}{TEXT}}

This technical subsection is dedicated to the study of the level curves $V(t_1,t_2)=v$, for $(t_1,t_2)\in\mathbf{P}$. We shall see that each of these level curves is a single connected arc, either entirely contained in $\mathbf{Q}^{\pm}$ or coinciding exactly with the hyperbola $\mathbf{H}$.

\begin{lem}\label{Lema:g_unimodal}
The function $g(t;v)$ given in \eqref{eq:L,g} is unimodal on $t>0$ for every $v>0$. Moreover the unique critical point $T_v$ of $g(t;v)$ is a global maximum point. 
\end{lem}
\begin{proof}
We have that  $\lim\limits_{t\to 0^+}g(t;v)=\lim\limits_{t\to +\infty}g(t;v)=0$ and thus, one critical point $T_v>0$ must exist. Moreover, we have:
\[\frac{\partial g}{\partial t}(t;v)=\frac{G_1(t;v)}{(t+v)^2},\quad G_1(t;v)=\frac{(t-1)(t+v)-(t^2+v)\log t}{(t-1)^2}.\]
Taking the derivative once more yields
\[\frac{\partial G_1}{\partial t}(t;v)=\frac{(t+v)}{t(t-1)^3}G_2(t),\quad G_2(t)=1-t^2+2 t\log t.\]

Since $G_2'$ attains its global maximum at $G_2'(1)=0$, we have that $G_2$ is strictly decreasing for $t>0$. Now, as $G_2(1)=0$, we have that $\operatorname{sign}(G_2(t))=\operatorname{sign}((1-t)^3)$. Computing $\lim\limits_{t\to 1}\partial_tG_1(t;v)=-(1+v)/3<0$, we conclude that $G_1$ is strictly decreasing and thus, $T_v$ must be unique and a global maximum point of $g(t;v)$.
\end{proof}

\begin{propo}\label{Propo:Cv}
Let $v>0$ and consider the level set $\mathcal{C}_v=\{(t_1,t_2)\in\mathbf{P}:V(t_1,t_2)=v\}$. Then there exists an analytic strictly decreasing bijection $\delta_v:(0,T_v)\to (T_v,+\infty)$, such that
\[\mathcal{C}_v=\{(t_1,\delta_v(t_1)):0<t_1<T_v\},\]
and 
\[\lim\limits_{t_1\to T_v^-}\delta_v(t_1)=T_v,\quad \lim\limits_{t_1\to 0^+}\delta_v(t_1)=+\infty.\]
Hence, $\mathcal{C}_v$ is path-connected. 
\end{propo}
\begin{proof}
As we have already stated in the proof of Proposition \ref{Propo:Reduction}, $(t_1,t_2)\in\mathcal{C}_v$ if and only if $g(t_1;v)=g(t_2;v)$. By Lemma \ref{Lema:g_unimodal}, $g(t;v)$ is unimodal for any given $v>0$ and since $\lim\limits_{t\to 0^+}g(t;v)=\lim\limits_{t\to +\infty}g(t;v)=0$, Lemma \ref{Lema:unimodal-involution} implies that the equation $g(t_1;v)=g(t_2;v)$ defines a unique involution $\delta_v$ on $(0,\infty)$. The desired bijection is the restriction $t_2=\delta_v(t_1)$ for $t_1\in(0,T_v)$. The limits in the statement follow from the aforementioned limits of $g(t;v)$ and the fact that $T_v$ is the fixed point of the involution. Hence, $\mathcal{C}_v$ is the graph of $\delta_v$ contained inside $\mathbf{P}$ and thus, is a path-connected set.
\end{proof}

The following result collects three equations that will be useful to determine the images of the maps $V(t_1,t_2)$ and $K(t_1,t_2)$ in each subdomain of the decomposition $\mathbf{P}=\mathbf{Q}^+\sqcup\mathbf{H}\sqcup\mathbf{Q}^-$. 

\begin{lem}%\label{lema:ident}
Let $(t_1,t_2)\in\mathbf{P}$. Then the following equations hold:
\begin{align}
V(t_1,t_2)-1 &=\frac{(t_1+1)(t_2+1)\left(g(t_1;1)-g(t_2;1)\right)}{L(t_2)-L(t_1)},\label{eq:lema:ident-V-1}\\
t_1t_2-V(t_1,t_2)&=\frac{t_1t_2\left(\Xi(t_2)-\Xi(t_1)\right)}{L(t_2)-L(t_1)},\label{eq:lema:ident-p-V}\\
K(t_1,t_2)+V(t_1,t_2)-1&=\frac{(\Xi(t_2)-\Xi(t_1))(t_1+V(t_1,t_2))(t_2+V(t_1,t_2))}{V(t_1,t_2)(L(t_2)-L(t_1))},\label{eq:lema:ident-K+V-1}
\end{align}
where $\Xi(t)=(t+1)L(t)/t=(t+1)\log t/(t-1)$. Moreover,
\[g(\exp(a);1)=\frac{a}{2\sinh(a)},\quad \Xi(\exp(a))=a\coth\left(\frac{a}{2}\right).\]
\end{lem}
\begin{proof}
%Recalling the definition of $V(t_1,t_2)$ in equation \eqref{eq:V,K}, the assertion that $V(t_1,t_2)>0$ follows immediately from the fact that $L(t_2)-L(t_1)>0$ and
%\[t_2L(t_1)-t_1L(t_2)=t_1t_2\left(\frac{L(t_1)}{t_1}-\frac{L(t_2)}{t_2}\right)=t_1t_2\left(L\left(\frac{1}{t_1}\right)-L\left(\frac{1}{t_2}\right)\right)>0.\]
The equations are directly obtained by elementary manipulations on the expressions of $V(t_1,t_2)$, $K(t_1,t_2)$, $L(t)$ and $\Xi(t)$.
\end{proof}

\begin{lem}\label{Lema:signs}
Let $(t_1,t_2)\in\mathcal{C}_v$. Then
\[{\rm sign}(v-1)={\rm sign}(t_1t_2-v)={\rm sign}(\log(t_1t_2)),\]
and thus, for $0<v<1$, $\mathcal{C}_v\subset\mathbf{Q}^-$, for $v>1$, $\mathcal{C}_v\subset\mathbf{Q}^+$ and $\mathcal{C}_1=\mathbf{H}$.
\end{lem}
\begin{proof}
Let $a_i=\log t_i$, so $a_1<a_2$. Since $|a_2|^2-|a_1|^2=a_2^2-a_1^2=(a_2+a_1)(a_2-a_1)$, we have that
\[{\rm sign}(|a_2|-|a_1|)={\rm sign}(a_1+a_2)={\rm sign}(\log(t_1t_2)).\]
Since $L(t_2)-L(t_1)>0$ and $V(t_1,t_2)=v$, by equation \eqref{eq:lema:ident-V-1}, we have that ${\rm sign}(v-1)={\rm sign}(g(t_1;1)-g(t_2;1))={\rm sign}(|a_2|-|a_1|)$ because $a\mapsto g(\exp(a);1)$ is even and strictly decreasing for $a>0$. Equation \eqref{eq:lema:ident-p-V} implies that ${\rm sign}(t_1t_2-v)={\rm sign}(\Xi(t_2)-\Xi(t_1))={\rm sign}(|a_2|-|a_1|)$ because $a\mapsto\Xi(\exp(a))$ is even and strictly increasing for $a>0$.

Hence, for $v>1$, every point $(t_1,t_2)\in\mathcal{C}_v$ must satisfy $t_1t_2>1$, which implies that $\mathcal{C}_v\subset\mathbf{Q}^+$. Analogously, for $0<v<1$, $\mathcal{C}_v\subset\mathbf{Q}^-$. Additionally, $V(t_1,t_2)=1$ if and only if $t_1t_2=1$, hence $\mathcal{C}_1=\mathbf{H}$.
\end{proof}

Using the above results, we can decompose each subdomain $\mathbf{Q}^-$, $\mathbf{H}$ and $\mathbf{Q}^+$ as a union of level sets of $V$. More precisely,
\[\mathbf{Q}^-=\bigcup_{0<v<1}\mathcal{C}_v,\quad \mathbf{H}=\mathcal{C}_1,\quad \mathbf{Q}^+=\bigcup_{v>1}\mathcal{C}_v,\]
each level curve being a connected arc. The next step is to prove that the function $K$ is injective along the level curves $\mathcal{C}_v$, when restricted to each one of the above subdomains. This will follow from the fact that the Jacobian determinant $\det D\Theta$ does not vanish in each $\mathbf{Q}^\pm$.

\subsection{The Jacobian of \texorpdfstring{$\Theta(t_1,t_2)$}{TEXT}}

Considering the map $\Theta$ defined in \eqref{eq:Theta}, one can directly compute its Jacobian determinant, which can be written, most conveniently, in the form:
\begin{equation}\label{eq:JTheta}
\det D\Theta=\frac{W(t_1,t_2)(t_2-t_1)^2}{(L(t_2)-L(t_1))t_1^2t_2^2}\left(t_1^2t_2^2 A_1(t_1,t_2)-V^2(t_1,t_2)A_2(t_1,t_2)\right),
\end{equation}
where 
\[W(t_1,t_2)=\frac{(t_1+V(t_1,t_2))(t_2+V(t_1,t_2))}{V(t_1,t_2)},\; A_1(t_1,t_2)=\frac{L'(t_1)+L'(t_2)}{(t_2-t_1)^2}-\frac{2(L(t_2)-L(t_1))}{(t_2-t_1)^3},\]
and
\[A_2(t_1,t_2)=\frac{(t_1+t_2)(L(t_2)-L(t_1))}{(t_2-t_1)^3}-\frac{(t_1L'(t_1)+t_2L'(t_2))}{(t_2-t_1)^2}.\]
For $(t_1,t_2)\in\mathbf{P}$, one can readily see that the factor $W(t_1,t_2)(t_2-t_1)^2/(t_1^2t_2^2(L(t_2)-L(t_1)))$ is strictly positive and thus, to determine the sign of the Jacobian determinant, we must consider the functions $A_1$ and $A_2$ in \eqref{eq:JTheta}. The next results provide a characterization of such functions as integrals of a rational function, which will allow us to reduce the determination of the sign of $\det D\Theta$ to a comparison between two moments of a single positive measure.

\begin{lem}\label{Lema:A1,A2}
For every $t>0$, $L(t)=\displaystyle\int_{0}^{\infty}\frac{t\,dx}{(1+x)(t+x)}$. Additionally, given $(t_1,t_2)\in\mathbf{P}$, define $Q(x)=(x+t_1)(x+t_2)$. Then
\begin{equation}\label{eq:A1,A2}
A_1(t_1,t_2)=\int_{0}^\infty\dfrac{x\,dx}{(1+x)Q^2(x)},\quad A_2(t_1,t_2)=\int_{0}^\infty\dfrac{x^2\,dx}{(1+x)Q^2(x)}.
\end{equation}
In particular, $A_1(t_1,t_2)>0$ and $A_2(t_1,t_2)>0$.
\end{lem}
\begin{proof}
For $t=1$, the first integral in the statement evaluates to $\int_{0}^\infty\frac{dx}{(1+x)^2}=1=L(1)$. For $t\neq 1$,
\[\int_0^\infty\frac{t\,dx}{(1+x)(t+x)}=\frac{t}{t-1}\int_0^\infty\left(\frac{1}{1+x}-\frac{1}{t+x}\right)\,dx=\frac{t}{t-1}\log t=L(t).\]
Now, since $L'(t)=\frac{1}{t-1}-\frac{\log t}{(t-1)^2}=\int_0^\infty\frac{x\,dx}{(1+x)(t+x)^2}$, we have that
\begin{align}
A_1(t_1,t_2)=&\frac{1}{(t_2-t_1)^2}\int_0^\infty\left(\frac{x}{(1+x)(t_1+x)^2}+\frac{x}{(1+x)(t_2+x)^2}\right)\,dx\nonumber\\
&\;-\frac{2}{(t_2-t_1)^3}\int_0^\infty\left(\frac{t_2}{(1+x)(t_2+x)}-\frac{t_1}{(1+x)(t_1+x)}\right)\,dx=\nonumber\\
=&\int_0^\infty\frac{1}{(t_2-t_1)^2}\left(\frac{x(t_2+x)^2+x(t_1+x)^2-2x\,Q(x)}{(1+x)Q^2(x)}\right)\,dx=\int_0^\infty\frac{x\,dx}{(1+x)Q^2(x)}.\nonumber
\end{align}
And finally,
\begin{align}
A_2(t_1,t_2)=&\frac{(t_1+t_2)}{(t_2-t_1)^3}\int_0^\infty\left(\frac{t_2}{(1+x)(t_2+x)}-\frac{t_1}{(1+x)(t_1+x)}\right)\,dx\nonumber\\
&\;-\frac{1}{(t_2-t_1)^2}\int_0^\infty\left(\frac{t_1x}{(1+x)(t_1+x)^2}+\frac{t_2x}{(1+x)(t_2+x)^2}\right)\,dx=\nonumber\\
=&\frac{1}{(t_2-t_1)^3}\int_0^\infty\frac{1}{(1+x)Q^2(x)}\left((t_1+t_2)(t_2(t_2+x)(t_1+x)^2-t_1(t_1+x)(t_2+x)^2)\right.\nonumber\\
    &\;\left.-(t_2-t_1)(t_1 x(t_2+x)^2+t_2 x(t_1+x)^2)\right)\,dx=\int_0^\infty\frac{x^2\,dx}{(1+x)Q^2(x)}.\nonumber
\end{align}
\end{proof}

We now introduce two auxiliary moment functions that will be used to provide a simpler expression for $\det D\Theta$.

\begin{lem}\label{Lema:B0,B1}
For $(t_1,t_2)\in\mathbf{P}$, define
\begin{equation}\label{eq:Bj}
B_j(t_1,t_2)=\int_0^\infty\frac{x^j\,dx}{(1+x)Q(x)},
\end{equation}
for $j=0,1$. Then
\[\frac{L(t_2)-L(t_1)}{t_2-t_1}=B_1(t_1,t_2),\quad \frac{t_2L(t_1)-t_1L(t_2)}{t_2-t_1}=t_1t_2B_0(t_1,t_2).\]
\end{lem}
\begin{proof}
The result follows from the integral formula for $L(t)$ given in Lemma \ref{Lema:A1,A2} and elementary manipulations such as in its proof.
\end{proof}

\begin{corol}\label{Corol:detDTheta}
Setting $\beta(t_1,t_2)=B_1(t_1,t_2)/B_0(t_1,t_2)$ and $\alpha(t_1,t_2)=A_2(t_1,t_2)/A_1(t_1,t_2)$, the Jacobian determinant of the map $\Theta$ is given by:
\begin{equation*}%\label{eq:JTheta2}
\det D\Theta=\frac{W(t_1,t_2)(t_2-t_1)^2V^2(t_1,t_2)}{(L(t_2)-L(t_1))t_1^2t_2^2}A_1(t_1,t_2)\left(\beta^2(t_1,t_2)-\alpha(t_1,t_2)\right).
\end{equation*}
\end{corol}
\begin{proof}
The result follows immediately from equation \eqref{eq:JTheta} and the fact that, by Lemma \ref{Lema:B0,B1}, $V(t_1,t_2)=t_1t_2 B_0(t_1,t_2)/B_1(t_1,t_2)$.
\end{proof}

We state a relevant result about a reflection property that will be useful in the proof of Theorem \ref{Teo:fibersTheta-complete}.

\begin{lem}\label{Lema:iota}
Consider the involution $\iota:\mathbf{P}\to\mathbf{P}$ given by $\iota(t_1,t_2)=(1/t_2,1/t_1)$. Then
\[B_j\circ\iota(t_1,t_2)=t_1 t_2 B_{1-j}(t_1,t_2),\quad A_\ell\circ\iota(t_1,t_2)=t_1^2t_2^2A_{3-\ell}(t_1,t_2),\]
for $j=0,1$ and $\ell=1,2$. Consequently
\[\alpha\circ\iota(t_1,t_2)=\frac{1}{\alpha(t_1,t_2)},\quad\beta\circ\iota(t_1,t_2)=\frac{1}{\beta(t_1,t_2)}.\quad \]
\end{lem}
\begin{proof}
The lemma follows easily from Lemma \ref{Lema:A1,A2} and direct computations. We exhibit here only one of them:
\[B_j\circ\iota(t_1,t_2)=\int_0^\infty\frac{x^j\,dx}{(1+x)(\frac{1}{t_2}+x)(\tfrac{1}{t_1}+x)}=\int_0^\infty\frac{t_1t_2w^{3-j}}{(1+w)Q(w)}\frac{dw}{w^2}=t_1t_2B_{1-j}(t_1,t_2).\]
\end{proof}

We now introduce a logarithmic change of variables that will allow us to determine the sign of $\beta^2(t_1,t_2)-\alpha(t_1,t_2)$. More precisely, it transforms each subdomain of the decomposition $\mathbf{P}=\mathbf{Q}^+\sqcup\mathbf{H}\sqcup\mathbf{Q}^-$ into a region determined by the sign of a single coordinate $\eta$. Additionally, it allows us to rewrite each $A_\ell$ and $B_j$ above as integrals with respect to the same measure, depending on $\eta$. The problem is then reduced to a one-dimensional analysis on the variable $\eta$.

More precisely, we will introduce the logarithmic change of variables $\Psi:\mathbf{P}\to\mathbb{R}\times\mathbb{R}_{+}$ given by $\Psi(t_1,t_2)=(\eta,\xi)$, where
\begin{equation}\label{eq:logchange}
\eta=\tfrac{1}{2}\log(t_1t_2),\quad \xi=\tfrac{1}{2}\log(t_2/t_1).
\end{equation}
Notice that $t_1=e^{\eta-\xi}$, $t_2=e^{\eta+\xi}$ and that $\mathbf{Q}^+=\Psi^{-1}(\{\eta>0\})$, $\mathbf{Q}^-=\Psi^{-1}(\{\eta<0\})$ and $\mathbf{H}=\Psi^{-1}(\{\eta=0\})$. Let us define the following auxiliary functions:
\begin{equation*}%\label{eq:vartheta,C,E,l}
\vartheta(u)=\frac{1}{1+e^{\eta+u}},\; C(u)=\frac{1}{\cosh u+\cosh\xi},\;E(u)=e^{u},\;\mathscr{L}(u)=\frac{C'(u)}{C(u)}=-\frac{\sinh u}{\cosh u+\cosh\xi}.
\end{equation*}
By substituting $x=\exp(\eta+u)$, $u\in\mathbb{R}$ into the formulas \eqref{eq:A1,A2} and \eqref{eq:Bj}, one obtains
\[B_j(t_1,t_2)=\frac{e^{(j-1)\eta}}{2}\int_{-\infty}^{+\infty} e^{ju}\vartheta(u)C(u)\,du,\quad A_\ell(t_1,t_2)=\frac{e^{(\ell-3)\eta}}{4}\int_{-\infty}^{+\infty}e^{(\ell-1)u}\vartheta(u)C^2(u)\,du.\]
Setting $d\mu=\vartheta(u)C(u)\,du$, we denote $\mom{\eta}{f}=\int_\mathbb{R}f\,d\mu$. It is important to note that the bracket operation $\mom{\eta}{\cdot}$ depends on both parameters $(\eta,\xi)$, although we highlight only its dependence on $\eta$ since that will play a substantial role in the next proofs. We have that
\[B_0(t_1,t_2)=\frac{e^{-\eta}}{2}\mom{\eta}{1},\quad B_1(t_1,t_2)=\frac{1}{2}\mom{\eta}{E},\]
\[ A_1(t_1,t_2)=\frac{e^{-2\eta}}{4}\mom{\eta}{C},\quad A_2(t_1,t_2)=\frac{e^{-\eta}}{4}\mom{\eta}{EC},\]
and
\begin{equation}\label{eq:beta,alpha}
\beta(t_1,t_2)=e^{\eta}\frac{\mom{\eta}{E}}{\mom{\eta}{1}},\quad \alpha(t_1,t_2)=e^{\eta}\frac{\mom{\eta}{EC}}{\mom{\eta}{C}}.
\end{equation}
It is important to remark that $1,E,C,EC,\mathscr{L}, C\mathscr{L}\in L^1(\mu)$. Indeed, since $\vartheta$ and $C$ are positive functions and $\vartheta(u)\leqslant\min\{1,e^{-\eta-u}\}$ and $C(u)\leqslant 2e^{-|u|}$, we have that
\[\mom{\eta}{1}\leqslant\int_{-\infty}^02e^u\,du+\int_0^{+\infty}2e^{-\eta-2u}\,du=2+e^{-\eta},\]
\[\mom{\eta}{E}\leqslant\int_{-\infty}^02e^{2u}\,du+\int_0^{+\infty}2e^{-\eta-u}\,du=1+2e^{-\eta}.\]
The bounds $C(u)<1$, $|\mathscr{L}(u)|<1$ imply that $\mom{\eta}{C}\leqslant\mom{\eta}{1}$ and $\mom{\eta}{EC}\leqslant \mom{\eta}{E}$ and $\mom{\eta}{|\mathscr{L}|}\leqslant\mom{\eta}{1}$ and $\mom{\eta}{ |C\mathscr{L}|}\leqslant\mom{\eta}{C}$. 

\begin{lem}[Differentiation formulas]
For a fixed $\xi>0$ and $f\in\{1,E,C,EC\}$, the function $\eta\mapsto \mom{\eta}{f}$ belongs to $C^1(\mathbb{R})$ and $\partial_\eta(\mom{\eta}{f})=-\int_\mathbb{R}\vartheta(u)(f C)'(u)\,du$. In particular,
\begin{align}\label{eq:diff_formulas1}
\partial_\eta(\mom{\eta}{1})&=-\mom{\eta}{\mathscr{L}},\quad \partial_\eta(\mom{\eta}{E})=-\mom{\eta}{E}-\mom{\eta}{E\mathscr{L}},\nonumber\\
\partial_\eta(\mom{\eta}{C})&=-2\mom{\eta}{C\mathscr{L}},\quad
  \partial_\eta(\mom{\eta}{E C})=-\mom{\eta}{EC}-2\mom{\eta}{E C\mathscr{L}},
\end{align}
Moreover, the following equations hold
\begin{equation}\label{eq:diff_formulas2}
\mom{\eta}{E\mathscr{L}}=-e^{-\eta}\mom{\eta}{\mathscr{L}},\quad \mom{\eta}{E C\mathscr{L}} =-e^{-\eta}\mom{\eta}{C\mathscr{L}},
\end{equation}
\begin{equation}\label{eq:diff_formulas3}
\mom{\eta}{\mathscr{L}}=\int_{\mathbb{R}}(-\vartheta'(u)C(u))\,du>0,\quad \mom{\eta}{C\mathscr{L}} =\frac{1}{2}\int_{\mathbb{R}}(-\vartheta'(u)C^2(u))\,du>0.
\end{equation}
\end{lem}
\begin{proof}
Note that $\partial_\eta\vartheta=\partial_u\vartheta=-\vartheta(1-\vartheta)$. Hence $|\partial_\eta\vartheta|\leqslant \vartheta$ and $\vartheta'(u)<0$. Thus, for $f\in\{1, E, C, EC\}$, $|\partial_\eta(\vartheta f C)|=|f C||\partial_\eta\vartheta|\leqslant \vartheta C|f|$. Since $\vartheta$ is decreasing in $\eta$, for $\eta$ in a compact interval $[\eta_0,\eta_1]$ the right-hand side is bounded by the integrable function $\vartheta\vert_{\eta=\eta_0}C|f|$, and the following differentiation under the integral sign is justified by dominated
convergence:
\[\partial_\eta(\mom{\eta}{f})=\int_\mathbb{R}f(u)C(u)\vartheta'(u)\,du=\left[f(u)C(u)\vartheta(u)\right]_{-\infty}^{+\infty}-\int_\mathbb{R}\vartheta(u)(fC)'(u)\,du.\]
Since $f(u)C(u)\vartheta(u)\to 0$ as $|u|\to\infty$, substituting each of the four possible functions into $f$ yields the equation \eqref{eq:diff_formulas1}.
Equation \eqref{eq:diff_formulas2} follows from the computation
\[\mom{\eta}{(e^{-\eta}+E)f}=\int_\mathbb{R}f(u)(e^u+e^{-\eta})\vartheta(u)C(u)\,du=e^{-\eta}\int_\mathbb{R}f(u)C(u)\,du.\]
For $f\in\{\mathscr{L},C\mathscr{L}\}$ the right-hand side of the above equation evaluates to $\left[e^{-\eta}\frac{C^j(u)}{j}\right]_{-\infty}^{\infty}=0$, with $j=1,2$. The final equation \eqref{eq:diff_formulas3} follows from the fact that $\vartheta'(u)<0$ and $C^j(u)\vartheta(u)\to 0$ as $|u|\to\infty$ for $j=1,2$, since
\[-\int_\mathbb{R}\vartheta'(u)\frac{C^j(u)}{j}\,du=\left[-\vartheta(u)\frac{C^j(u)}{j}\right]_{-\infty}^{+\infty}+\int_\mathbb{R}\vartheta(u)C^{j-1}(u)C'(u)\,du=\int_\mathbb{R}\mathscr{L}(u)\vartheta(u)C^j(u)\,du.\]
\end{proof}

\begin{corol}\label{Corol:partial_eta,beta,alpha}
For a fixed $\xi>0$, consider $\beta$ and $\alpha$ given in \eqref{eq:beta,alpha} as functions of $\eta$. Then $\alpha,\beta\in C^1(\mathbb{R})$, and
\begin{equation}\label{eq:diff_formulas_alpha,beta}
\partial_\eta\beta=(1+\beta)\frac{\mom{\eta}{\mathscr{L}}}{\mom{\eta}{1}},\quad \partial_\eta\alpha=2(1+\alpha)\frac{\mom{\eta}{C\mathscr{L}}}{\mom{\eta}{C}}.
\end{equation}
%where
%\begin{equation}\label{eq:m,n}
%\mathbf{m}=\frac{\mom{\eta}{\mathscr{L}}}{\mom{\eta}{1}}>0,\quad \mathbf{n}=\frac{\mom{\eta}{C\mathscr{L}}}{\mom{\eta}{C}}>0.
%\end{equation}
\end{corol}
\begin{proof}
%The fact that $\mathbf{m},\mathbf{n}>0$ is immediate from equation \eqref{eq:diff_formulas3}. 
By the differentiation formulas \eqref{eq:diff_formulas1} and \eqref{eq:diff_formulas2}, we have
\[\partial_\eta(\log\beta)
   =1+\frac{\partial_\eta\mom{\eta}{E}}{\mom{\eta}{E}}
     -\frac{\partial_\eta\mom{\eta}{1}}{\mom{\eta}{1}}
   =\frac{\mom{\eta}{\mathscr{L}}}{\mom{\eta}{1}}
     \Bigl(\frac{e^{-\eta}\mom{\eta}{1}}{\mom{\eta}{E}}+1\Bigr)=\frac{\mom{\eta}{\mathscr{L}}}{\mom{\eta}{1}}\Bigl(\frac{1}{\beta}+1\Bigr).\]
Analogously,
\[ \partial_\eta\log\alpha
   =1+\frac{\partial_\eta\mom{\eta}{E C}}{\mom{\eta}{E C}}-\frac{\partial_\eta\mom{\eta}{C}}{\mom{\eta}{C}}=2\Bigl(e^{-\eta}\frac{\mom{\eta}{C\mathscr{L}}}{\mom{\eta}{EC}}+\frac{\mom{\eta}{C\mathscr{L}}}{\mom{\eta}{C}}\Bigr)=2\frac{\mom{\eta}{C\mathscr{L}}}{\mom{\eta}{C}}\Bigl(\frac{1}{\alpha}+1\Bigr).
\]
\end{proof}

The equations \eqref{eq:diff_formulas_alpha,beta} in the previous result form a decoupled 1-parametric family of differential equations which determine $\beta$ and $\alpha$ completely from their initial value at $\eta=0$. One can see that, since the coefficients are positive, both $\beta$ and $\alpha$ grow with $\eta$. The sign of $\beta^2(t_1,t_2)-\alpha(t_1,t_2)$ is then determined by which of the two grows faster.

The next step in our chain of arguments is to prove that 
\[\frac{\mom{\eta}{C\mathscr{L}}}{\mom{\eta}{C}}<\frac{\mom{\eta}{\mathscr{L}}}{\mom{\eta}{1}}.\]
In order to do so, we will need to recall the classical \emph{Korkine's identity} (see, for instance \cite{Kumar04,MPF}).

\begin{lem}[Korkine's identity] Let $m:\mathbb R\to[0,\infty)$ be measurable with $\int_{\mathbb R} m(x)\,dx<\infty$. Let $f_1,f_2:\mathbb R\to\mathbb R$ be measurable with $mf_1,\ mf_2,\ mf_1f_2\in L^1(\mathbb R)$. Then
\begin{align}\label{eq:Korkine}
\int_{\mathbb R} m(x)\,dx\int_{\mathbb R} f_1(x)f_2(x)m(x)\,dx-\int_{\mathbb R} f_1(x)m(x)\,dx\int_{\mathbb R} f_2(x)m(x)\,dx
=\nonumber\\
=\frac{1}{2}\iint_{\mathbb{R}^2}\big(f_1(x)-f_1(y)\big)\big(f_2(x)-f_2(y)\big)m(x)m(y)\,dx\,dy ,
\end{align}
where the right-hand side is an absolutely convergent integral over $\mathbb R^2$ with respect to $x$ and $y$.
\end{lem}

We remark that the expression in the left-hand side of \eqref{eq:Korkine} is known in the literature as the \emph{weighted Chebyshev functional} (see \cite{Bezziou19}, for instance). 

\begin{propo}\label{Propo:n<m}
For every $(\eta,\xi)\in\mathbb{R}\times\mathbb{R}_{+}$, we have that 
\[\frac{\mom{\eta}{C\mathscr{L}}}{\mom{\eta}{C}}<\frac{\mom{\eta}{\mathscr{L}}}{\mom{\eta}{1}}.\]
\end{propo}
\begin{proof}
The inequality in the statement is equivalent to $\mom{\eta}{C\mathscr{L}}\mom{\eta}{1}-\mom{\eta}{C}\mom{\eta}{\mathscr{L}}<0$. Since $1,C,\mathscr{L},C\mathscr{L}\in L^1(\mu)$, applying Korkine's identity \eqref{eq:Korkine}, the result holds once we prove that the following integral is negative:
\[I=\iint_{\mathbb{R}^2}Z(x,y)m(x)m(y)\,dx\,dy,\]
where $m(u)=\vartheta(u)C(u)$ and $Z(x,y)=(C(x)-C(y))(\mathscr{L}(x)-\mathscr{L}(y))$. As $\mathscr{L}$ is an odd function and $C$ is even, we see that $Z(-x,-y)=-Z(x,y)$. Moreover, $\mathscr{L}'(u)=-\frac{(1+\cosh u\cosh\xi)}{(\cosh u+\cosh\xi)^2}<0$ implies that $\mathscr{L}$ is strictly decreasing on $\mathbb{R}$. Additionally, since $C'(u)=-\frac{\sinh u}{(\cosh u+\cosh\xi)^2}<0$, $C(u)$ is also strictly decreasing for $u>0$. Thus, for $x+y>0$, with $x\neq y$, ${\rm sign}(x-y)={\rm sign}(|x|-|y|)$ and we obtain
\[{\rm sign}(C(x)-C(y))={\rm sign}(y-x)={\rm sign}(\mathscr{L}(x)-\mathscr{L}(y)),\]
hence $Z(x,y)>0$. The antisymmetry $Z(-x,-y)=-Z(x,y)$ implies that for $x\neq y$, ${\rm sign}(Z(x,y))={\rm sign}(x+y)$, while $Z(x,x)=0$, by definition.

We claim that ${\rm sign}(m(x)m(y)-m(-x)m(-y))=-{\rm sign}(x+y)$. In fact, consider the function $\chi(u)=\log(m(-u))-\log(m(u))$. Since $m(u)=\vartheta(u)C(u)$, the fact that $C$ is even implies that $\chi(u)=\log(1+e^{u+\eta})-\log(1+e^{-u+\eta})$ and thus $\chi'(u)=\frac{e^{u+\eta}}{1+e^{u+\eta}}+\frac{e^{\eta-u}}{1+e^{\eta-u}}>0$. Therefore, $\chi$ is a strictly increasing and odd function. For $x+y>0$, one obtains
\[0<\chi(x)-\chi(-y)=\chi(x)+\chi(y)=\log\left(\frac{m(-x)m(-y)}{m(x)m(y)}\right),\]
which is equivalent to $m(x)m(y)-m(-x)m(-y)<0$. The claim is then proved by symmetry.
Now, by the change of variables $(x,y)\mapsto(-x,-y)$, we have that 
\[I=-\iint_{\mathbb{R}^2}Z(x,y)m(-x)m(-y)\,dx\,dy,\]
which yields
\[2I=\iint_{\mathbb{R}^2}Z(x,y)\left(m(x)m(y)-m(-x)m(-y)\right)\,dx\,dy.\]
The integrand is non-positive on $\mathbb{R}^2$ and negative in the set $\{x+y\neq 0\}\cap\{x\neq y\}$ which is a set of full measure. Therefore $I<0$.
\end{proof}

We are now in the position to determine the sign of the Jacobian determinant $\det D\Theta$.

\begin{figure}[!ht]
    \centering
    \vspace*{0.1cm}
    \begin{overpic}[width=0.75\linewidth]{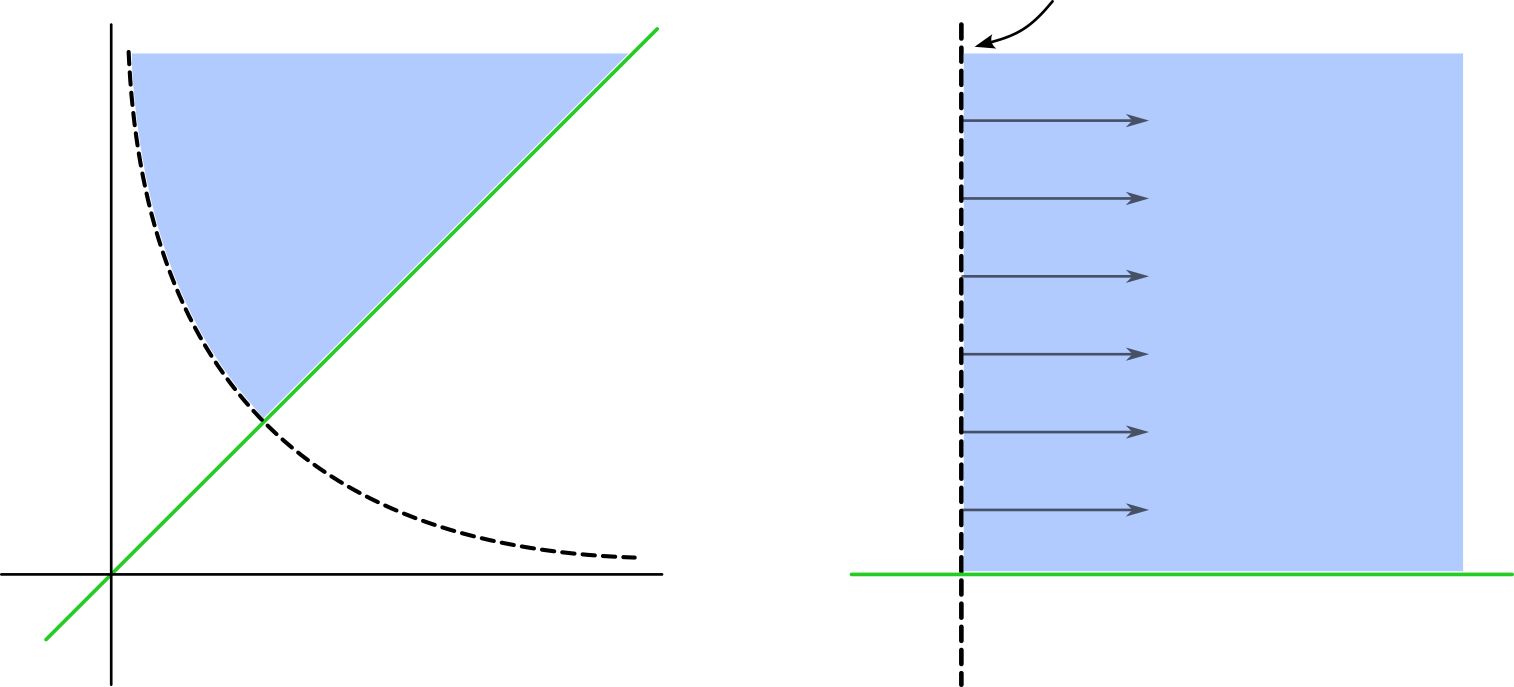}
    \put(15,32){$\mathbf{Q}^+$}
    \put(5,45){$t_2$}
    \put(44,6){$t_1$}
    \put(101,6){$\eta$}
    \put(62,45){$\xi$}
    \put(70,46){$\det D\Theta=0$}
    \put(72,23.5){$\det D\Theta>0$}
    \end{overpic}
    \caption{The schematics of the proof of Theorem \ref{Teo:JacobianTheta}: The Jacobian determinant is zero at the set $t_1t_2=1$ which in the coordinates $(\eta,\xi)$ is given by $\eta=0$. Integrating the functions $\beta$ and $\alpha$ along the curves $\xi\equiv\xi_0$ allows us to determine the sign of the Jacobian determinant on the quadrant $\eta>0$ which corresponds to $\mathbf{Q}^+$ in the coordinates $(t_1,t_2)$.}
    \label{Fig:Propagation}
\end{figure}

\begin{teo}\label{Teo:JacobianTheta}
For the map $\Theta$ given by \eqref{eq:Theta} defined on the domain $\mathbf{P}$, we have that
\[{\rm sign}(\det D\Theta)={\rm sign}(\log(t_1t_2)).\]
\end{teo}
\begin{proof}
Consider the change of variables \eqref{eq:logchange}. Thus, $\log(t_1t_2)=2\eta$ and the theorem will be proven once we consider the three cases $\eta=0$, $\eta>0$ and $\eta<0$, which correspond, respectively, to $(t_1,t_2)\in\mathbf{H}$, $(t_1,t_2)\in\mathbf{Q}^+$ and $(t_1,t_2)\in\mathbf{Q}^-$.

The case $\eta=0$, i.e. $(t_1,t_2)\in\mathbf{H}$ follows immediately from Lemma \ref{Lema:iota}, since for $t_1t_2=1$, we have $\iota(t_1,t_2)=(t_1,t_2)$, and therefore
$\alpha(t_1,t_2)=\beta(t_1,t_2)=1$. Hence, by Corollary \ref{Corol:detDTheta}, we have that $\det D\Theta(t_1,t_2) = 0$ for $(t_1,t_2)\in\mathbf{H}$.\medskip

%At $\eta=0$, we have $\vartheta(-u)=(1+e^{-u})^{-1}=e^u\vartheta(u)$. Thus, for any given $f\in L^1(\mathbb{R})$, the following equation holds
%\[\mom{\eta}{f(-\;\cdot)}=\int_\mathbb{R}f(-u)\vartheta(u)C(u)du=\int_\mathbb{R}f(u)\vartheta(-u)C(u)du=\int_\mathbb{R}e^uf(u)\vartheta(u)C(u)du=\mom{\eta}{Ef}.\]
%In particular, for $f\in\{1,C\}$, since both functions are even, we have that $\mom{\eta}{1}=\mom{\eta}{E}$ and $\mom{\eta}{C}=\mom{\eta}{E C}$, which imply that $\beta=\alpha=1$. Hence, by Corollary \ref{Corol:detDTheta}, we have that $\det D\Theta(t_1,t_2) = 0$ for $(t_1,t_2)\in\mathbf{H}$.

Let us fix $\xi>0$ and consider the curve $\eta\mapsto (t_1,t_2)=(e^{\eta-\xi},e^{\eta+\xi})$ inside $\mathbf{P}$ which crosses $\mathbf{H}$ exactly at $\eta=0$ (see Figure \ref{Fig:Propagation}). By Corollary \ref{Corol:partial_eta,beta,alpha},
\[\partial_\eta\log(1+\beta)=\frac{\mom{\eta}{\mathscr{L}}}{\mom{\eta}{1}},\quad \partial_\eta\log(1+\alpha)=2\frac{\mom{\eta}{C\mathscr{L}}}{\mom{\eta}{C}}.\]
Integrating the above equations with initial condition $\beta=\alpha=1$ at $\eta=0$ yields
\[\beta=2e^\mathscr{B}-1,\quad \alpha=2e^{2\mathscr{A}}-1,\]
where
\[\mathscr{A}=\int_0^\eta\frac{\mom{s}{C\mathscr{L}}}{\mom{s}{C}}\,ds,\quad\mathscr{B}=\int_0^\eta\frac{\mom{s}{\mathscr{L}}}{\mom{s}{1}}\,ds.\]
Now we consider the case $\eta>0$. By Proposition \ref{Propo:n<m} and Corollary \ref{Corol:partial_eta,beta,alpha}, we have that $0<\frac{\mom{\eta}{C\mathscr{L}}}{\mom{\eta}{C}}<\frac{\mom{\eta}{\mathscr{L}}}{\mom{\eta}{1}}$, which implies that $0<\mathscr{A}<\mathscr{B}$. Hence,
\[\beta^2-\alpha=(2e^{\mathscr{B}}-1)^2-2e^{2\mathscr{A}}+1=2(e^\mathscr{B}-1)^2+2(e^{2\mathscr{B}}-e^{2\mathscr{A}})>0.\]
Again, by Corollary \ref{Corol:detDTheta}, $\det D\Theta(t_1,t_2)>0$ for $(t_1,t_2)\in\mathbf{Q}^+$.

The case $\eta<0$ follows directly by applying the involution $\iota(t_1,t_2)=(1/t_2,1/t_1)$ which is taken to $(-\eta,\xi)$ via $\Psi$, with $-\eta>0$. The previous case implies that $(\beta\circ\iota)^2-\alpha\circ\iota>0$ and by Lemma \ref{Lema:iota},
\[\beta^2-\alpha=\frac{1}{(\beta\circ\iota)^2}-\frac{1}{\alpha\circ\iota}=\frac{\alpha\circ\iota-(\beta\circ\iota)^2}{(\alpha\circ\iota)\,(\beta\circ\iota)^2}<0.\]
Thus, Corollary \ref{Corol:detDTheta} implies that $\det D\Theta(t_1,t_2)<0$ for $(t_1,t_2)\in\mathbf{Q}^-$.
\end{proof}

\subsection{The image of the map \texorpdfstring{$\Theta$}{TEXT}}

In this subsection we will determine the image of the map $\Theta:\mathbf{P}\to\mathbb{R}^2$ on each of its subdomains. In order to do so, we first need to define the critical curve $k=k_{\rm crit}(v)$ which will be crucial to write $\Theta(\mathbf{P})$ as a disjoint union. Recall that, by Lemma \ref{Lema:g_unimodal} and Proposition \ref{Propo:Cv}, for every $v>0$, there exists a unique $T_v$ which is a critical point of $g(\cdot;v)$ and for which the level curve $\mathcal{C}_v$ is given by the graph of an analytic strictly decreasing function $\delta_v:(0,T_v)\to(T_v,+\infty)$. We define the \emph{critical curve} by the equation $k=k_{\rm crit}(v)$, where
\begin{equation*}%\label{eq:k_c}
k_{\rm crit}(v)=\frac{T_v^2}{v}+2T_v-\frac{v(2T_v+v)}{T_v^2}. 
\end{equation*}

Intuitively, the critical curve has the following significance: one can see, by direct computation, that for $(t_1,t_2)\in\mathcal{C}_v$, the equation $K(t_1,t_2)=k$ is equivalent to $q(t_1;v,k)=q(t_2;v,k)$, where
\[q(t;v,k)=\frac{t^3+vkt+v^2}{t(t+v)}.\]
This function is unimodal on $t>0$, with $q(t;v,k)\to+\infty$ as $t\to0^+$ and $t\to+\infty$. Thus, by Lemma \ref{Lema:unimodal-involution}, both equations in \eqref{eq:JetSepConditions_simplified} are defined by involutions of $(0,+\infty)$. The critical curve is precisely the set of parameters for which the fixed points of both these involutions coincide. Indeed, $\partial_t q(T_v;v,k)=0$ if and only if $k=k_{\rm crit}(v)$.

\begin{lem}\label{Lema:k_cinverte}
For any $v>0$, $T_{1/v}=1/T_v$ and $k_{\rm crit}(1/v)=-k_{\rm crit}(v)/v$. In particular $k_{\rm crit}(1)=0$.
\end{lem}
\begin{proof}
We have that $g(t;1/v)=\dfrac{t\log t}{(t-1)(t+1/v)}=\dfrac{v\,t^{-1}\log t}{(1-t^{-1})(v+t^{-1})}=v\,g(t^{-1};v)$. Hence
\[\frac{\partial g}{\partial t}(T_v^{-1};1/v)=-v\,{T_v^2}\,\frac{\partial g}{\partial t}(T_v;v)=0,\]
and thus $T_{1/v}=1/T_v$. Substituting into $k_{\rm crit}(1/v)$, we obtain
\[k_{\rm crit}(1/v)=\frac{1}{v}\left(\frac{v^2}{T_v^2}+\frac{2v}{T_v}-2T_v-\frac{T_v^2}{v}\right)=-\frac{k_{\rm crit}(v)}{v}.\]
\end{proof}

\begin{lem}\label{Lema:K_Cv}
Consider the restriction $K\vert_{\mathcal{C}_v}:\mathcal{C}_v\to\mathbb{R}$ which is parametrized by $K(t_1,\delta_v(t_1))$. Then, the following holds:
\begin{itemize}
    \item[(i)] $\lim\limits_{t_1\to T_v^-}K\vert_{\mathcal{C}_v}=k_{\rm crit}(v)$;
    \item[(ii)] For $v>1$, $\lim\limits_{t_1\to 0^+}K\vert_{\mathcal{C}_v}=+\infty$.
    %\item[(iii)] For $0<v<1$, $\lim\limits_{t_1\to 0^+}K\vert_{\mathcal{C}_v}=-\infty$;
\end{itemize}
\end{lem}
\begin{proof}
For every $(t_1,t_2)\in \mathcal{C}_v$, we have that $K(t_1,t_2)=\dfrac{t_1t_2}{v}+t_1+t_2-\dfrac{v(t_1+t_2+v)}{t_1t_2}$, with $t_2=\delta_v(t_1)$. Proposition \ref{Propo:Cv} implies that $t_2\to T_v$ as $t_1\to T_v^-$ which proves item (i).

To prove (ii), we consider equation \eqref{eq:lema:ident-K+V-1}. Since $v>1$, by Lemma \ref{Lema:signs}, we have that $t_1t_2>v$ and thus
\[\Xi(t_2)-\Xi(t_1)>\log t_2-\frac{t_1+1}{t_1-1}\log t_1=\log(t_1t_2)-\frac{2t_1\log t_1}{t_1-1}>\log v-2L(t_1).\]
Moreover, $(t_1+v)(t_2+v)>t_2v$, $0<L(t_2)-L(t_1)<L(t_2)$ and, by Proposition \ref{Propo:Cv}, $t_2\to+\infty$ as $t_1\to 0^{+}$. Since $L(t_1)\to 0$ as $t_1\to0^+$, we have $\log v-2L(t_1)>0$ for $t_1$ small enough, and then
\[K\vert_{\mathcal{C}_v}>1-v+\left(\log v-2L(t_1)\right)\frac{t_2}{L(t_2)}\to+\infty.\]
\end{proof}

The previous lemma allows us to see the critical curve at the limit of $\Theta(t_1,t_2)$ as $(t_1,t_2)\to (t,t)$ (see Figure \ref{Fig:ImageTheta}). We can now describe completely the map $\Theta$. We finish the section by proving the following result which is a more complete formulation of Theorem \ref{Teo:fibersTheta}.\medskip

\begin{teo}\label{Teo:fibersTheta-complete}
Consider the map $\Theta:\mathbf{D}\to\mathbb{R}^2$ given by \eqref{eq:Theta} and the decomposition $\mathbf{P}=\mathbf{Q}^+\sqcup\mathbf{H}\sqcup\mathbf{Q}^-$. Then, the following holds:
\begin{itemize}
    \item[(a)] $\mathbf{H}=\mathcal{C}_1$ and $\Theta\vert_\mathbf{H}\equiv (1,0)$;
    \item[(b)] $\Theta$ restricts to a diffeomorphism on $\mathbf{Q}^+$:
    \[\Theta:\mathbf{Q}^+\to \{(v,k): v>1, k>k_{\rm crit}(v)\};\]
    \item[(c)] $\Theta$ restricts to a diffeomorphism on $\mathbf{Q}^-$:
    \[\Theta:\mathbf{Q}^-\to \{(v,k): 0<v<1, k<k_{\rm crit}(v)\}.\]
\end{itemize}
Moreover for every
$(v,k)\in\Theta(\mathbf Q^+)\cup\Theta(\mathbf Q^-)$, the fiber
$\Theta^{-1}(v,k)$ consists of exactly two points, namely
$(t_1,t_2)$ and $(t_2,t_1)$ for a unique $(t_1,t_2)\in\mathbf P$.
\end{teo}
\begin{proof}
(a): Lemma \ref{Lema:signs} yields $\mathbf{H}=\mathcal{C}_1$. Now, let $(t,1/t)$ parametrize $\mathbf{H}$. A direct computation yields $V(t,1/t)=1$ and thus $K(t,1/t)=1+t+1/t-(t+1/t+1)=0$.

In order to prove (b) and (c), we will first prove that $\Theta\vert_{\mathbf{Q}^\pm}$ is injective. In fact, let $z_1,z_2\in\mathbf{Q}^{\pm}$ such that $\Theta(z_1)=\Theta(z_2)$. Then, $z_1,z_2\in \mathcal{C}_v$ for a common $v\neq 1$. By Proposition \ref{Propo:Cv}, we can set $\gamma(t)=(t,\delta_v(t))$ as the parametrization of the level curve $\mathcal{C}_v$, and then there exist $s_1,s_2\in(0,T_v)$ such that $z_i=(s_i,\delta_v(s_i))$. Thus,
\[(V\circ\gamma)'(t)=\nabla V(\gamma(t))\cdot\gamma'(t)\equiv 0,\;\text{for all }t\in(0,T_v).\]
Consider $K\circ\gamma$. Since $\mathcal{C}_v$ is path-connected, if $s_1\neq s_2$, then Rolle's Theorem would imply the existence of $t^\ast\in(s_1,s_2)$ such that $(K\circ\gamma)'(t^\ast)=0$. In that case,
\[(K\circ\gamma)'(t^\ast)=\nabla K(\gamma(t^\ast))\cdot\gamma'(t^\ast)=0,\]
which would imply that the gradients of $V$ and $K$ are parallel at $\gamma(t^\ast)$ and then $\det D\Theta (\gamma(t^\ast))=\nabla V(\gamma(t^\ast))\wedge\nabla K(\gamma(t^\ast))=0$. By Theorem \ref{Teo:JacobianTheta}, this is a contradiction and thus $z_1=z_2$. In fact, Theorem \ref{Teo:JacobianTheta} implies that $(K\circ\gamma)'(s)\neq 0$ for every $s\in(0,T_v)$. Thus, $K\vert_{\mathcal{C}_v}$ is monotone.

Now, for $(t_1,t_2)\in\mathbf{Q}^+$, Lemma \ref{Lema:signs} implies $\Theta(t_1,t_2)=(v,k)$ with $v>1$. Since $K\vert_{\mathcal{C}_v}$ is monotone we have, by Lemma \ref{Lema:K_Cv}, $k\in(k_{\rm crit}(v),+\infty)$, which proves (b).

To prove (c), we consider the involutions $\iota:\mathbf{P}\to\mathbf{P}$ and $\tau:\mathbb{R}_{+}\times\mathbb{R}\to\mathbb{R}_{+}\times\mathbb{R}$ given by $\iota(t_1,t_2)=(1/t_2,1/t_1)$ and $\tau(v,k)=(1/v,-k/v)$. Then, 
\[V\circ\iota(t_1,t_2)=\frac{\tfrac{1}{t_2}L(1/t_1)-\tfrac{1}{t_1}L(1/t_2)}{L(1/t_2)-L(1/t_1)}=\frac{L(t_1)-L(t_2)}{t_1L(t_2)-t_2L(t_1)}=\frac{1}{V(t_1,t_2)},\]
\[K\circ\iota(t_1,t_2)=\frac{V(t_1,t_2)}{t_1t_2}+\frac{t_1+t_2}{t_1t_2}-\frac{t_1t_2}{V(t_1,t_2)}\left(\frac{t_1+t_2}{t_1t_2}+\frac{1}{V(t_1,t_2)}\right)=-\frac{K(t_1,t_2)}{V(t_1,t_2)},\]
that is, $\Theta\circ\iota=\tau\circ\Theta$. Then, since $\iota(\mathbf{Q}^-)=\mathbf{Q}^+$, we have that $\Theta\vert_{\mathbf{Q}^-}=\tau\circ\Theta\circ\iota$ is a diffeomorphism by item (b) above. Its image is the set of points $(w,l)=\tau(v,k)$, with $v>1$ and $k>k_{\rm crit}(v)$, and hence $w=1/v<1$ and $l=-k/v<-k_{\rm crit}(v)/v=k_{\rm crit}(w)$ by Lemma \ref{Lema:k_cinverte}. 
The last assertion follows from the fact that $\Theta(t_1,t_2)=\Theta(t_2,t_1)$. This completes the proof.
\end{proof}

\section*{Acknowledgments}

\noindent L. Q. Arakaki is supported by Funda\c c\~ao de Amparo \`a Pesquisa do Estado de S\~ao Paulo (grant numbers 2024/06926-7 and 2025/08206-4). D. Panazzolo is partially supported by the ANR Project NonSper ANR-23-CE40-0028. For the purpose of Open Access, the authors have applied a CC-BY public copyright licence to any Author Accepted Manuscript (AAM) version arising from this submission.

\section*{Declaration of generative AI and AI-assisted technologies in the manuscript preparation process}

\noindent 

During the preparation of this work the authors used Claude (Opus 5 and Opus 5.5) to assist in the proofs of the technical lemmas, namely the rewriting of the quantities $A_i$, $B_j$ as integrals of rational functions and the transport arguments in the proof of Theorem \ref{Teo:JacobianTheta}. After using this tool/service, the authors reviewed and edited the content as needed and take full responsibility for the content of the published article.

\section*{Conflict of interest}

\noindent The authors declare that they have no conflict of interest.

\section*{Data availability}

\noindent Data sharing is not applicable to this article as no new data were created or analyzed in this study.

\bibliographystyle{siam}
\bibliography{Ref.bib}

\end{document}